\documentclass[12pt]{article}
\usepackage{amsmath}
\usepackage{amsfonts}
\usepackage{epsfig}
\newcommand{\CP}{\mathbb{CP}}
\newcommand{\C}{\mathbb{C}}
\newcommand{\Z}{\mathbb{Z}}

\newcommand{\p}{\pi_1}

\newcommand{\vp}{\varphi}
\newcommand{\ri}{\rightarrow}
\newcommand{\e}{\ell}
\newcommand{\tp}{\tilde{p}}
\newcommand{\usr}{\underset\sim\rightarrow}
\newcommand{\quadf}{\quad\quad\quad\quad\quad\quad}
\begin{document}
\newtheorem{corollary}{Corollary}[section]
\newtheorem{lemma}{Lemma}[section]
\title{THE REGENERATION OF A 5-POINT}
\author{MICHAEL FRIEDMAN AND MINA TEICHER}
\date{\today}
\maketitle

\begin{abstract}
The braid monodromy factorization of the branch curve of a surface
of general type is known to be an invariant that completely
determines the diffeomorphism type of the surface (see \cite{KuTe}).
Calculating this factorization is of high technical complexity;
computing the braid monodromy factorization of branch curves of
surfaces uncovers new facts and invariants of the surfaces. Since
finding the branch curve of a surface is very difficult, we
degenerate the surface into a union of planes. Thus, we can find the
braid monodromy of the branch curve of the degenerated surface,
which is a union of lines. The regeneration of the singularities of
the branch curve, studied locally, leads us to find the global braid
monodromy factorization of the branch curve of the original surface.
So far, only the regeneration of the BMF of 3,4 and 6-point (a
singular point which is the intersection of 3 / 4 / 6 planes; see
\cite{MoTe2},\cite{MoTe4}) were done. In this paper, we fill the gap
and find the braid monodromy of the regeneration of a 5-point. This
is of great importance to the understanding of the BMT (braid
monodromy type) of surfaces \cite{KuTe}.

 This braid monodromy will
be used to find the global braid monodromy factorization of
different surfaces; in particular - the monodromy of the branch
curve of the Hirzebruch surface $F_{2,(2,2)}$.

\end{abstract}

\tableofcontents

\section[Introduction]{Introduction}

Let $X \subset \CP^N$ be a smooth algebraic surface of degree $n$.
One may obtain information on $X$ by considering it as a branched
cover of another surface. If the base surface is $\CP^2$ and if the
map $X \rightarrow \CP^2$ is a generic projection, then the branch
locus is a plane curve $\overline S \subset \CP^2$ which is, in
general, singular. If the projection is generic, the singularities
are nodes and cusps. Let $ S \subset \C^2 \subset \CP^2$ denote a
generic affine portion of $\overline S$. A general problem is the
study of the fundamental groups of the complement of the branch
curve: $\p(\C^2 - S)$ and $\p(\CP^2 - \overline S)$.

It has been proven that these fundamental groups (derived from braid
monodromy factorizations) are invariants that distinguish between
diffeomorphic surfaces (see \cite{KuTe}); that is, if two surfaces
have equivalent braid monodromy factorizations (and thus isomorphic
fundamental groups), then they are diffeomorphic. However, the
converse is not true, that is, the diffeomorphism type does not
determine the equivalence class of the factorization. In \cite{KhKu}
a pair of diffeomorphic surfaces was constructed such that the braid
monodromy factorizations are not equivalent.

The above fundamental groups cannot be found directly, since finding
the branch curve explicitely is very difficult. Therefore, one has
to degenerate the surface $X$ into a union of planes, where in this
case, the branch curve is easy to find -- it is an arrangment of
lines. It is known, by the Zariski-Van Kampen Theorem that the braid
monodromy factorization (BMF ; see Section 2 for its definition) of
the branch curve determines the desired fundamental groups. Note
that the BMF of any curve is given by a product of the local BMF in
the neighborhood of the singular points of the branch curve. Thus,
the BMF of any line arrangement can be found explicitly (see
\cite{MoTe1}). By applying the regeneration techniques on the
singularities of the arrangement of lines, one can find the BMF of
the original branch curve.

So it is very important to find out what are the local BMF that are
obtained from regenerating different line arrangements (or line and
conic arrangement).

Till now, most of the arrangements that include one line and one
conic (or two lines) were treated. In \cite{MoTe2} the BMFs of the
regeneration of a tangent / node / branch point are given. However,
for more complicated arrangements, only a few results are known. The
BMF of the regeneration of a standard 3-point (that is, a singular
point which is the intersection of 3 planes), 4-point and 6-point
are presented in \cite{MoTe4}. Figure I.1 is a depiction of the
arrangement of planes (that correspond to the regions delimited by
the edges of the diagram) and edges (that correspond to lines of
intersection between two planes).
\begin{center}
\epsfig{file=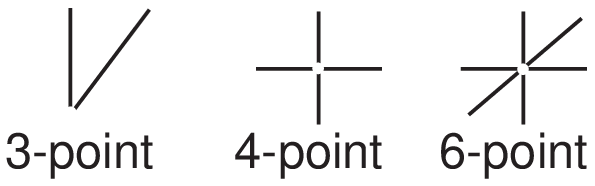}\\
(figure I.1)
\end{center}
Notice that in the cases in which the line arrangements include more
than two lines, the order of the regeneration effects how the
factorization will look.

In this article we compute two important braid monodromy
factorizations, which were not known till now -- the BMF of a
5-point, and a general formula of a certain type of $(k+1)$-point
where $k \geq 2$.

This article is organized as follows: In Section 2 we give the main
definitions (BMF and regeneration techniques), and then we compute
the BMF of the two main cases that were mentioned above. Section 3
shows the importance of these factorizations by introducing an
example which uses one of them.

\textbf{Acknowledgment}: The authors wish to thank Prof. Eugenii
Shustin for his help and for fruitful discussions.

\section[5-point regeneration]{5-point regeneration}
This section introduces the main result of the article -- the local
braid monodromy factorization induced from the regeneration of a
neighborhood of a 5-point; that is, a point which is the
intersection of 5 planes. Our result deals with two specific cases
for which this situation can appear, though there are other
configurations of 5-planes passing through a point. Note that we
actually consider this point to be a singular point of the branch
curve of a degenerated surface, when considering its generic
projection to $\CP^2$. But first we need to recall a few
definitions, related to the braid monodromy factorization and to the
regeneration techniques.

\subsection[Preliminaries]{Preliminaries: BMF and regeneration techniques}

Computing the braid monodromy is the main tool to compute
fundamental groups of complements of curves. In this subsection we
define the braid monodromy.\\


Let $D$ be a closed disk in $ \mathbb{R}^2,$ \ $K\subset Int(D),$
$K$ finite, $n= \#K$. Recall that the braid group $B_n[D,K]$ can be
defined as the group of all equivalent diffeomorphisms $\beta$ of
$D$ such that $\beta(K) = K\,,\, \beta |_{\partial D} =
\text{Id}\left|_{\partial D}\right.$. \\

\textit{Definition}:\ \underbar{$H(\sigma)$, half-twist defined by
$\sigma$}

Let $a,b\in K,$ and let $\sigma$ be a smooth simple path in $Int(D)$
connecting $a$ with $b$ \ s.t. $\sigma\cap K=\{a,b\}.$ Choose a
small regular neighborhood $U$ of $\sigma$ contained in $Int(D),$
s.t. $U\cap K=\{a,b\}$. Denote by $H(\sigma)$ the diffeomorphism of
$D$ which switches $a$ and $b$ by a counterclockwise 180 degree
rotation and is the identity on $D\setminus U$\,. Thus it defines an
element of $B_n[D,K],$ called {\it the half-twist defined by
$\sigma$ }.\\
%
%
%

Assume that all of the points of $K$ are on the $X$-axis (when
considering $D$ in $\mathbb{R}^2$). In this situation, if $a,b \in
K$, and $z_{a,b}$ is a path that connects them, then we denote it by
$Z_{a,b} = H(z_{a,b})$. If $z_{a,b}$ is a path that goes below the
$x$-axis, then we denote it by $\underline Z_{a,b}$, or just
$Z_{a,b}$. If $z_{a,b}$ is a path that goes above the $x$-axis, then
we denote it by $\overline Z_{a,b}$. See \cite{MoTe2}, Section 2 for
additional
notations.\\

 \textit{Definition}: \underline{The braid monodromy w.r.t.
 $S,\pi,u$}

 Let $S$ be a curve, $S\subseteq \C^2$
 Let $\pi: S\to\C^1$ be defined by
 $\pi(x,y)=x.$
We denote $\deg\pi$ by $m.$ Let $N=\{x\in\C^1\bigm| \#\pi^{-1}(x)<
m\}.$
   Take $u\notin N,$ s.t.  $\Re(x)\ll u$ \ $\forall x\in N.$
Let  $ \C^1_u=\{(u,y)\}.$  There is a  natural defined homomorphism
$$\pi_1(\C^1-N,u)\xrightarrow{\vp} B_m[\C_u^1,\C_u^1\cap S]$$ which
is called {\it the braid monodromy w.r.t.} $S,\pi,u,$ where $B_m$ is
the braid group. We sometimes denote $\vp$ by $\vp_u$. In fact,
denoting by $E$ -- a big disk in $\C^1$ s.t. $E \supset N$, we can
also take the path in $E\setminus N$ not to be a loop, but just a
non-self-intersecting path; this induces a diffeomorphism between
the models $(D,K)$ at the two ends of the considered path, where $D$
is a big disk in $\C^1_u$, and
$K = \C_u^1\cap S \subset D$.\\

\textit{Definition}:  $\underline{\psi_T, \ \text{Lefschetz
diffeomorphism induced by a path} \ T }$

Let  $T$ be a path in $E\setminus N$ connecting $x_0$ with $x_1$,
$T:[0, 1]\ri E\setminus N$. There exists a continuous family of
diffeomorphisms $\psi_{(t)}: D\ri D,\ t\in[0,1],$ such that
$\psi_{(0)}=Id$, $\psi_{(t)}(K(x_0))=K(T(t)) $ for all $t\in[0,1]$,
and  $\psi_{(t)}(y)= y$ for all $y\in\p D$. For emphasis we write
$\psi_{(t)}:(D,K(x_0))\ri(D,K(T(t))$. Lefschetz diffeomorphism
induced by a path $T$ is the diffeomorphism
$$\psi_T= \psi_{(1)}: (D,K(x_0))\usr (D,K(x_1)).$$
 Since $ \psi_{(t)} \left( K(x_{0})\right) = K(T(t))$ for all
$t\in [0,1]$, we have a family of canonical isomorphisms
$$\psi_{(t)}^{\nu}: B_p\left[ D, K(x_{0})\right] \usr B_p\left[
D, K({T(t)})\right], \ \quad \text{for all} \, \, t\in[0,1].$$\\
see the following figure for illustration of the above definitions:
\begin{center}
\epsfig{file=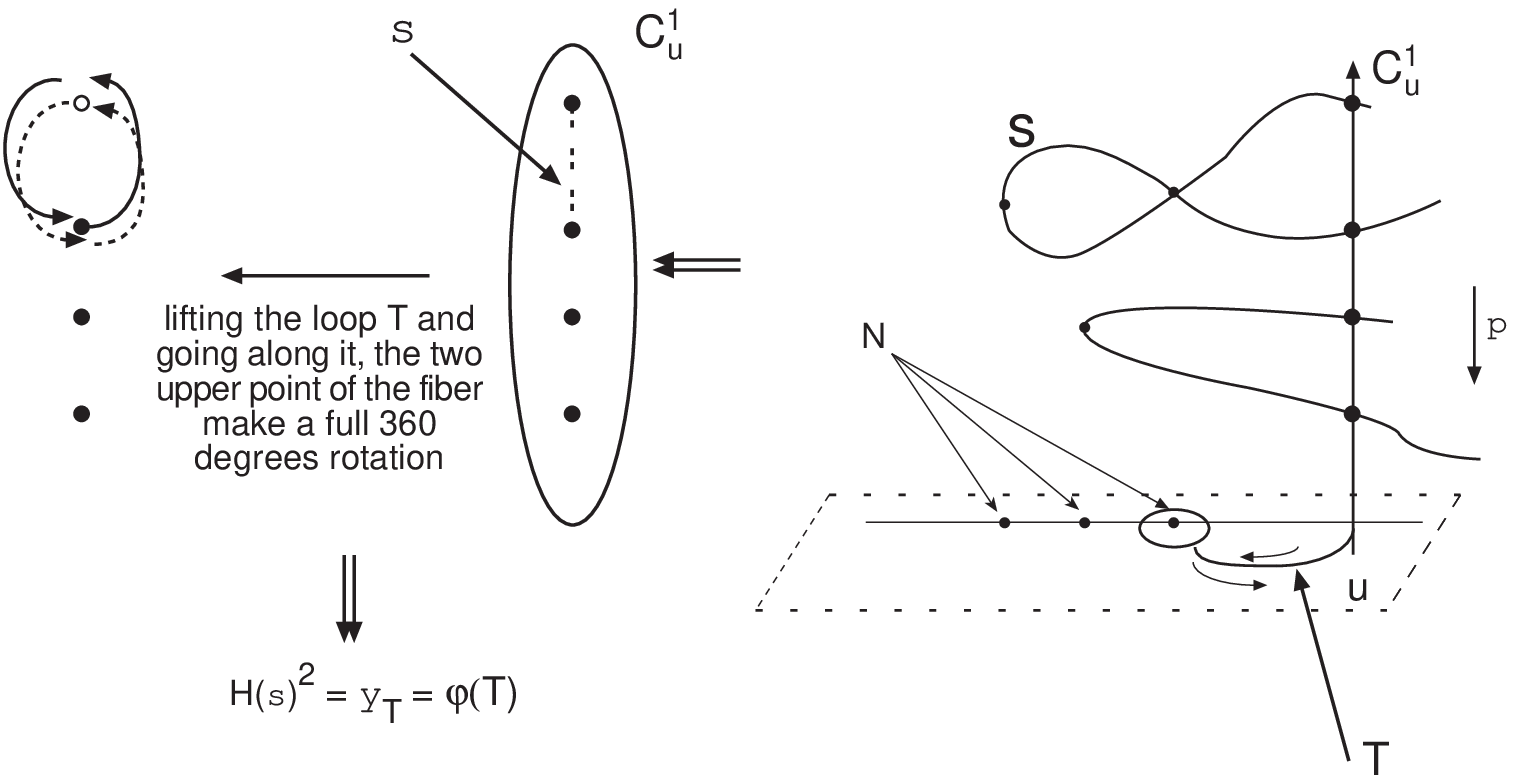}\\
\end{center}

We recall Artin's theorem on the presentation of the Dehn twist of
the braid group as a product of braid monodromy elements of a
geometric-base (a base of $\p = \p(\C^1 - N, u)$ with certain
properties; see \cite{MoTe1} for definitions).\\
 \textbf{Theorem}: Let $S$ be a curve transversal to
the line in infinity, and $\vp$ is a braid monodromy of $S , \vp:\p
\rightarrow B_m$. Let {$\delta_i$} be a geometric (free) base
(g-base) of $\p,$ and $ \Delta^2$ is the generator of Center($B_m$).
Then:
$$\Delta^2 = \prod\vp(\delta_i).$$ This product is also defined as
the \textsl{braid monodromy factorization} (BMF) related to a curve $S$.\\


So in order to find out what is the braid monodromy factorization of
$\Delta_p^2$, we have to find out what are $\vp (\delta_i),\,\forall
i$. We refer the reader to the definition of a \textit{skeleton}
(see \cite{MoTe2}) $\lambda_{x_j}, x_j \in N$, which is a model of a
set of paths connecting points in the fiber, s.t. all those points
coincide when approaching $A_j=$($x_j,y_j$)$\in S$, when we approach
this point from the right. To describe this situation in greater
detail, for $x_j \in N$, let $x_j' = x_j + \alpha$. So the skeleton
in $x_j$ is defined as system of paths connecting the points in
$K(x_j') \cap D(A_j,\varepsilon)$ when $0 < \alpha \ll \varepsilon
\ll 1$, $D(A_j,\varepsilon)$ is a disk centered
in $A_j$ with radius $\varepsilon$.\\

For a given skeleton, we denote by
$\Delta\langle\lambda_{x_j}\rangle$ the braid by rotates by 180
degrees counterclockwise a small neighborhood of the given skeleton.
Note that of $\lambda_{x_j}$ is a single path, then
$\Delta\langle\lambda_{x_j}\rangle = H(\lambda_{x_j})$.

We also refer the reader to the definition of $\delta_{x_0}$, for
$x_0 \in N$ (see \cite{MoTe2}), which describes the Lefschetz
diffeomorphism induced by a path going below $x_0$, for different
types of singular points (tangent, node, branch; for example, when
going below a node, a half-twist of the skeleton occurs. When going
below a tangent point, a full-twist occurs).

We define, for $x_0 \in N$, the following number: $\varepsilon_{x_0}
= 1,2,4$ when ($x_0, y_0$) is a branch / node / tangent point
(respectively). So we have the following statement (see
\cite{MoTe2}, prop. 1.5):

 Let $\gamma_j$ be a path below the real
line from $x_j$ to $u$, s.t. $\ell(\gamma_j)=\delta_j$. So -
$$\vp_u(\delta_j) = \vp(\delta_j) =
\Delta<(\lambda_{x_j})(\prod\limits_{m=j-1}^{1}\delta_{x_m})>^{\varepsilon_{x_j}}.$$
When denoting $\xi_{x_j} =
(\lambda_{x_j})\Bigg(\prod\limits_{m=j-1}^{1}\delta_{x_m}\Bigg)$ we
get --
$$\vp(\delta_j)
= \Delta\langle(\xi_{x_j})\rangle^{\varepsilon_{x_j}}.$$ Note that
the last formula gives an algorithm to compute the wanted
factorization.

For a detailed explanation of the braid monodromy, see \cite{MoTe1}.\\

 We recall now the regeneration
methods.

 The regeneration methods are actually, locally, the
reverse process of the degeneration method. When regenerating a
singular configuration consisting of lines and conics, the final
stage in the regeneration process involves doubling each line, so
that each point of $K$ corresponding to a line labelled $i$ is
replaced by a pair of points, labelled $i$ and $i'$. The purpose of
the regeneration rules is to explain how the braid monodromy behaves
when lines are doubled in this manner. We denote by $Z_{i,j} =
H(z_{i,j})$ where $z_{i,j}$ is a path connecting points in $K$.

 The rules are (see \cite{MoTe4}, pg. 336-7):
\begin{enumerate}
\item \textbf{First regeneration rule}: The regeneration of a branch point of hyperbola:\\
A factor of the braid monodromy of the form $Z_{ij}$ is replaced in
the regeneration by $Z_{i'j}\cdot
\overset{(j)}{\underline{Z}}_{ij'}$
\item \textbf{Second regeneration rule}:The regeneration of a node:\\
A factor of the form $Z^2_{ij}$ is replaced by a factorized
expression $Z^2_{ii',j} := Z^2_{i'j}\cdot Z^2_{ij}$ , $Z^2_{i,jj'}
:= Z^2_{ij'}\cdot Z^2_{ij}$ or by $Z^2_{ii',jj'} := Z^2_{i'j'}\cdot
Z^2_{ij'}Z^2_{i'j}\cdot Z^2_{ij}$.
\item \textbf{Third regeneration rule}:The regeneration of a tangent point:\\
A factor of the form $Z^4_{ij}$ in the braid monodromy factorized
expression is replaced by $Z^3_{i,jj'} := (Z^3_{ij})^{Z_{jj'}}\cdot
(Z^3_{ij}) \cdot (Z^3_{ij})^{Z^{-1}_{jj'}}$.
\end{enumerate}
As a result, we get a factorized expression, which, by \cite{KuTe},
determines the diffeomorphism type of our surface, and, by
\cite{VK}, determines $\pi_1(\CP^2 -\overline S)$.
%

\subsection[The first case]{The first case}

In this subsection we will look at the case where, locally, we have
5 planes corresponding to the angular sectors of the figure,
intersecting each other along a line whenever they have a common
edge. The lines $L_i,\,1\leq i \leq 5$ are numerated as following:
\begin{center}
\epsfig{file=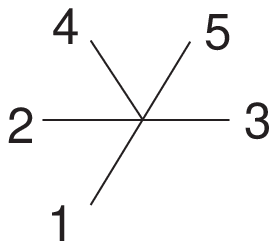}\\
(figure 1)
\end{center}
The lines are numerated in a way such that it describes the
respective positions of the points where they intersect
$\mathbb{C}^1_u$. We know that a 5-point of this sort can be
generated during the process of a degeneration of a surface into a
union of planes. Thus, when we examine the local braid monodromy
factorization of this 5-point, before degenerating, we know it is
$\Delta^2\langle1,5\rangle$. So by knowing what the regeneration
process will do to this factorization, we will know part of the
relations which are in the local fundamental group of $\C^2$ (or
$\CP^2$) minus the branch curve in this local neighborhood. Note
that in the regeneration process, line 4 is regenerated first, then
lines 2 and 3 and then lines 1 and 5.

In order to compute the desired factorization, we need a few
corollaries. The first is cited from \cite{MoTe4}, and deals with
the result of the regeneration process under certain conditions.
\begin{corollary}
Let $V$ be a projective algebraic surface, $D'$ -- a curve in $V$.
Let $f:V \rightarrow \mathbb{CP}^2$ be a generic projection. Let
$S\subseteq \mathbb{CP}^2, S'\subseteq V$ be the corresponding
branch / ramification curve of $f$. Assume $S'$ intersects $D'$ in
$\alpha'$. Let $D=f(D'), \alpha = f(\alpha')$. Assume that there
exist neighborhoods of $\alpha$ and $\alpha'$ s.t. $f|_{S'}$ and
$f|_{D'}$ are isomorphic. Then $D$ is tangent to $S$ at $\alpha$.
\end{corollary}
\textbf{Proof:} see \cite{MoTe4}.\\

The second corollary deals with the computation of a few braids,
which are induced from loops going around a complex intersection of
a conic and a line. We need this lemaa, since this situation appears
during the regeneration process. So consider the following model.

Let $C=\{(y^2-x)(y+x+1)=0\},\,\pi_1,\pi_2:C\to\mathbb{C}
,\,\pi_1(x,y)\mapsto{x},\,\pi_2(x,y)\mapsto{y}$. Denote by $p_1,
p_2$ the points of intersection of $y^2=x$ and $y=-x-1$. So --
$x_{p_1}=-\frac{1}{2}+\frac{\sqrt{3}i}{2},
x_{p_2}=-\frac{1}{2}-\frac{\sqrt{3}i}{2}$. Denote --
$x_0=-\frac{1}{4},\,
A=\pi_2(\pi_1^{-1}(x_0))=\{\pm{\frac{1}{2}i},-\frac{3}{4}\},\,
x_1=-\frac{3}{4},\,
A'=\pi_2(\pi_1^{-1}(x_1))=\{\pm\frac{\sqrt{3}}{2}i,-{\frac{1}{4}}\}$
(see figure 2).
\begin{center}
\epsfig{file=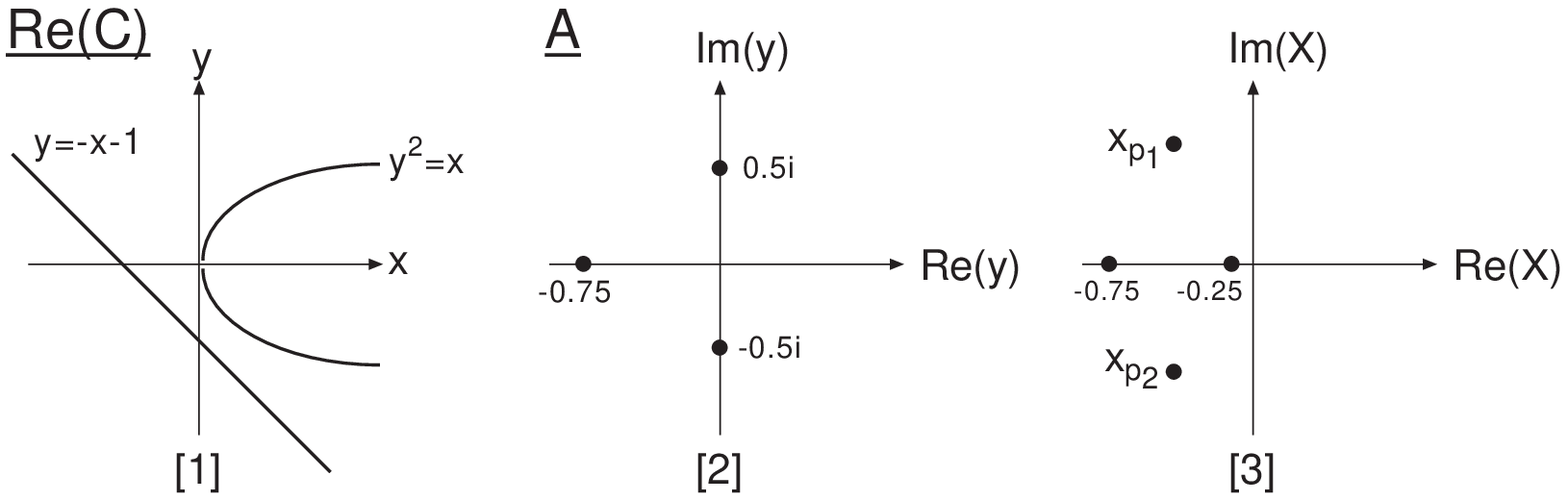}\\
(figure 2)
\end{center}
\textbf{Remark}:\:Note that $A'$ and $A$ (which are on the $Y$-axis)
are equivalent in the sense that if $a\in A, a' \in A'$ and $\Re(
a), \Re (a')\neq0$ (or $\Im(a), \Im(a') > 0$ or $\Im (a), \Im (a') <
0)$), then $a$ and $a'$ come from the same component of $C$.  Let
$D$ be a disk on the Y-axis s.t. $A,A'\subset D$. Thus we can define
a continuous diffeomorphism $\beta:D\rightarrow D$ s.t.
$\beta(\frac{-3}{4})=\frac{-1}{4},\beta(\frac{\pm i}{2})=\frac{\pm
\sqrt{3}i}{2}, \beta(\partial D)= \partial D, \beta(\Im(Y) \cap D) =
\Im(Y) \cap D , \beta(\Re(Y) \cap D) = \Re(Y) \cap D$, and outside a
small neighborhood of the $\Im(Y)$-- and $\Re(Y)$--axis, $\beta =
Id$.

Denote:\\
$\sigma_1$ -- the segment connecting $-\frac{3}{4}$ and
$\frac{i}{2}$ in $A$;\\$\sigma_2$ -- the segment connecting
$-\frac{3}{4}$ and $\frac{-i}{2}$ in $A$; See figure
3.[1]\\$\sigma'_1$ -- the segment connecting $-\frac{1}{4}$ and
$\frac{\sqrt{3}i}{2}$ in $A'$;\\$\sigma'_2$ -- the segment
connecting $-\frac{1}{4}$ and $-\frac{\sqrt{3}i}{2}$ in $A'$.
\\ Let $l_1(t), l_2(t), (0\leq t\leq1)$ be two loops starting (and
ending) at $x_0$, s.t. $l_i$ is around $x_{p_i}$. Lifting $l_i
(i=1,2)$ to $C$ and projecting it to D, we get a motion:
$(D,A)\rightarrow (D,A)$, which induces a braid; let $l_3(t), (0\leq
t\leq1)$ be a curve starting at $x_1$, ending at $x_0$ and
surrounding $x_{p_2}$ from below (see figures 3.[2], 3.[3]).
\begin{center}
\epsfig{file=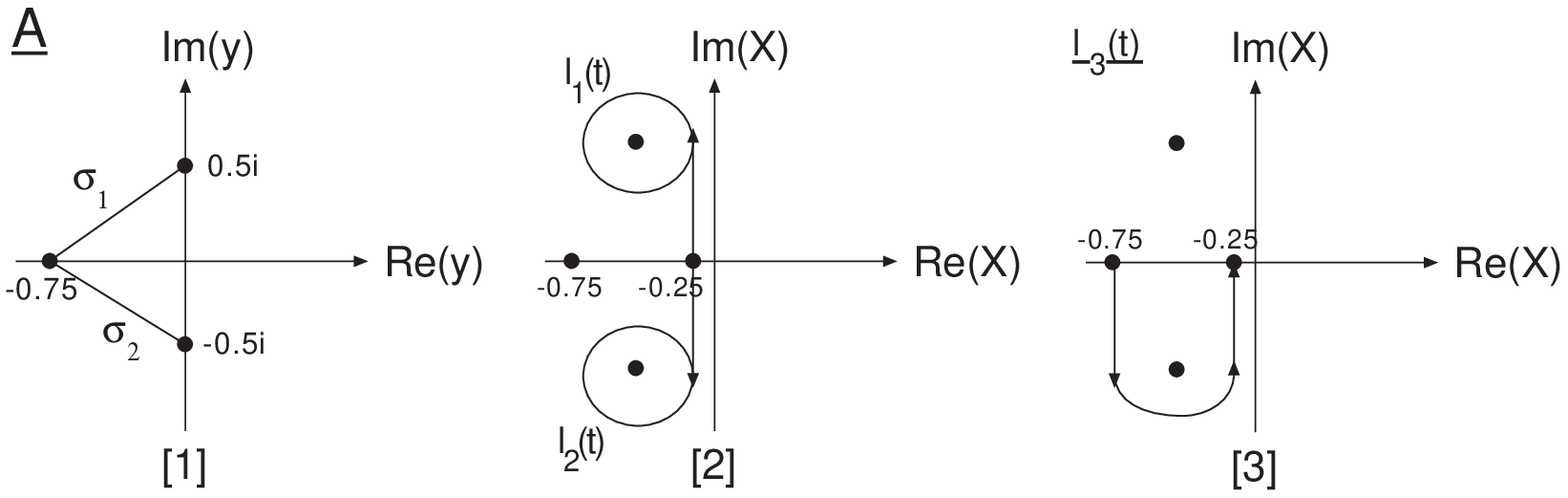}\\
(figure 3)
\end{center}
As above, we get a motion $(D,A)\rightarrow (D,A')$. Denote:\\
$\psi_{l_1},\psi_{l_2}$ -- the motions induced from $l_1, l_2$ (that
is, the Lefschetz isomorphisms induced by the paths; we omit
the superscript, as mentioned in the introduction).\\
$\psi_{l_3}$ -- the motion induced from $l_3$ (after composing
$\beta$ on the resulting disk).
\begin{corollary}
$\psi_{l_1} = H(\sigma_2)^2,\psi_{l_2}= H(\sigma_1)^2,\psi_{l_3}=
H(\sigma'_1)^2$, where $H(\sigma)$ is the halftwist induced from the
path $\sigma$.
\end{corollary}
\textbf{Proof:}\: For $\psi_{l_1},\psi_{l_2}$ we can look at a small
neighborhood of $p_1 \,(p_2)$. Since in this small neighborhood the
intersection of the branch of the conic, which $p_1\, (p_2)$ lies
on, and the line can be treated as the intersection of 2 lines, it
is obvious that when we perform a full loop around $x_{p_1}$ or
around $x_{p_2}$, the induced motion of the above points (points
$0.5i$ and $-0.75$ or points $-0.5i$ and $-0.75$) is a rotation of
360 degrees, and the induced braid is $H(\sigma)^2$ where $\sigma$
is the path connecting the points ($\sigma_2$ or $\sigma_1$; see
figure 3.[1]). Indeed, the line is $y=-x-1$, so when $x$ moves along
$l_1$, the corresponding value of $y$ is in the lower half-plane;
both end-points of $\sigma_2$ approach $-\frac12-\frac{\sqrt3}{2}i$
as $x$ approaches $x_{p_1}$; And similarly for $l_2$.

We shall now compute $\psi_{l_3}$. Observe that the union of the
straight line from $x_0$ to $x_1$ with the arc $l_3$ is a closed
loop from $x_0$ to itself, homotopic to $l_2$. Moving $x$ along the
real axis from $x_0$ to $x_1$ induces the diffeomorphism
$\beta:(D,A) \to (D,A')$ introduced earlier, so (up to isotopy)
$\psi_{l_3}\circ\beta=\psi_{l_2}$, which gives $\psi_{l_3}$. \quad
\quadf\quadf\quadf\quadf $\square$

So we now look at a point $v$, which is the intersection of 5 lines
(see figure 1), which are (part of a) branch curve of a degenerated
surface. Since we consider that this branch curve is a result of a
degeneration process, we can apply Corollary 2.1  when we are trying
to find out what will happen when we first regenerate line 4.

Denote $V = \bigcup\limits_{i=1}^5 L_i$. Define
$\underset{\{k\}}{\underline{Z}_{i\,j}} =
H(\underset{\{k\}}{\underline{z}_{i\,j}})$, where
$\underset{\{k\}}{\underline{z}_{i\,j}}$ is the path from point $i$
to point $j$, when the part of the path which is between $i$ and $j$
is below the $X$-axis, and it surrounds point $k$ from the left (if
$k<i$) or from the right (if $j<k$). For example, see the following
figures:\\
\begin{center}$\underset{\{7\}}{\underline{z}_{2\,5}}$
\epsfig{file=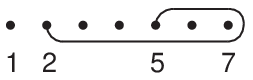}\\
$\underset{\{1\}}{\underline{z}_{2\,5}}$
\epsfig{file=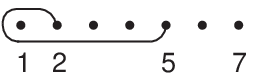}\\
\end{center}

\begin{corollary} After regenerating $V$ in a small neighborhood
$U$ of $v, \,L_4$ turns into a conic $Q_4$, s.t. $Q_4$ is tangent to
$L_2$ and $L_5$. Denote the resulting branch curve, after this
regeneration, by $\tilde{V}$. Thus, the singularities of $T =
\tilde{V} \cap U$ are as in the figure below:
\begin{center}
\epsfig{file=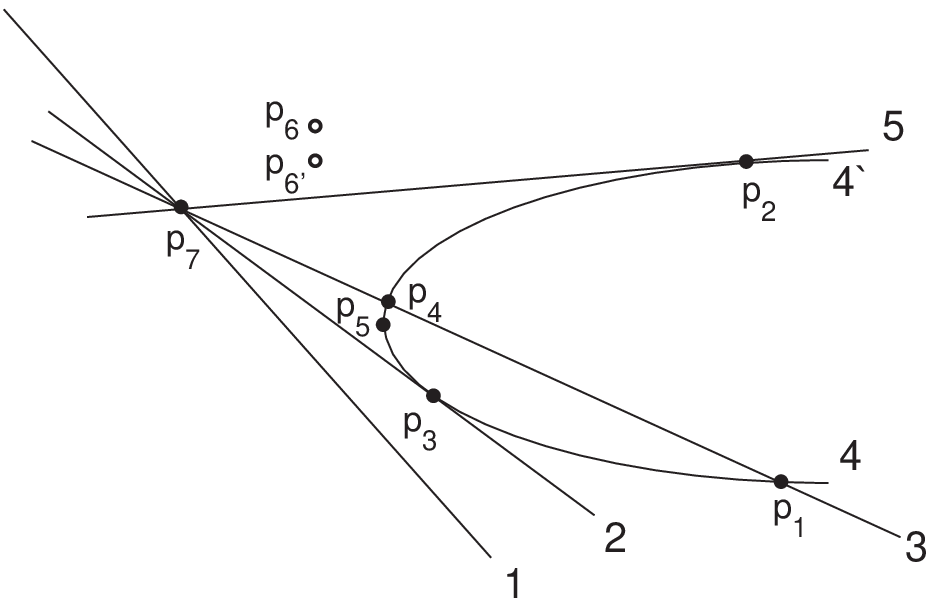}\\
\emph{(figure 4)}
\end{center}
Then the local braid monodromy of the above configuration is \\
$$\tilde{\varphi} = Z^2_{3\,4}\,Z^4_{4'\,5}Z^4_{2\,4}\,
\underset{\{5\}}{\overset{(4)}{\underline{Z}^2_{3\,4'}}}\,\hat{\hat{Z}}_{4\,4'}
\,\underset{\{5\}}{\underline{Z}^2_{1\,4'}}\,
\underline{Z}^2_{1\,4}\,
(\Delta^2\!<\!1,2,3,5\!>)^{Z^{-2}_{4\,5}},$$ where
$\hat{\hat{Z}}_{4\,4'} = H(\hat{\hat{z}}_{4,4'})
\,(\underset{\{5\}}{\overset{(4)}{\underline{Z}^2_{3\,4}}}=
H(\underset{\{5\}}{\overset{(4)}{\underline{z}^2_{3\,4}}}))$ is the
half-twist corresponding to the following path
:\\$\hat{\hat{z}}_{4,4'}$ (figure 5.[1]),\,
$\underset{\{5\}}{\overset{(4)}{\underline{z}^2_{3\,4}}}$ (figure
5.[2]):
\begin{center}
\epsfig{file=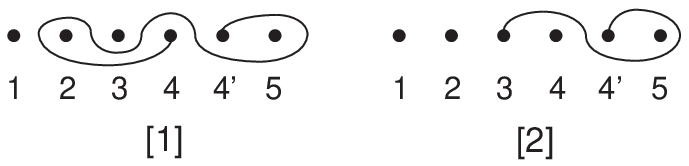}\\
\emph{(figure 5)}
\end{center}
\end{corollary}

\textbf{Proof:} Let $\{p_j\}_{j=1}^7 \cup \{p_{6'}\}$ be the
singular points of a small neighborhood (that is - $U$) of $v$ (see
figure 4) with respect to $\pi_1$ (the projection to the
X-axis) as follows: \\
$p_2, p_3$ -- the tangent points of $Q_4$ and $L_2, L_5$.\\
$\{p_1,p_4\},\{p_6,p_{6'}\}$ -- the intersection points of $Q_4$
with
$L_3,L_1$.\\
$p_5$ -- the branch point of $Q_4$.\\
$p_7$ -- the intersection point of $\{L_i\}_{i=1,2,3,5}$. \\
Let E (resp. D) be a closed disk on the $X$-axis (resp. $Y$-axis).
Let $N = \{x(p_j) = x_j \:|\: 1\leq j\leq 7 $\:or\:$ j = 6' \}$,
s.t. $N\subset E -
\partial E$. Let $M$ be a real point on the $x$-axis, s.t. $x_j \ll M, \forall x_j \in N,
 1\leq j\leq 7 $\:or\:$ j = 6' $. There is a $g$-base
$\ell (\gamma_j)_{j=1}^7 \cup \ell (\gamma_{6'})$ of $\pi_1 (E -
N,u)$, s.t. each path $\gamma_j$ is below the real line and the
values of $\varphi_M$ with respect to this base and $E \times D$ are
the ones given in the proposition. We look for
$\varphi_M(\ell(\gamma_j))$ for $j=1,\ldots,7$ or $ j = 6'$. Choose
a $g$-base $\ell\{\gamma_j\}_{j=1}^{7}\cup \ell (\gamma_{6'})$
as above and put all the data in the following table:\\

%
\begin{center}
\begin{tabular}{cccc}
$j$ & $\lambda_{x_j}$ & $\varepsilon_{x_j}$ & $\delta_{x_j}$
\\[0.5ex]\hline
1 & $<3,4>$ & 2 & $\Delta<3,4>$\\
2 & $<4',5>$ & 4 & $\Delta^2<4',5>$\\
3 & $<2,3>$ & 4 & $\Delta^2<2,3>$\\
4 & $<4,4'>$ & 2 & $\Delta<4,4'>$\\
5 & $<3,4>$ & 1 & $\Delta^{\frac{1}{2}}_{IR}<2>$\\
6,6' & $<1,3>,<1,4>$ & 2 & $\Delta^2<1,3>$\\
7 & $<1,2,4',5>$ & 2 & $ - $
\end{tabular}
\end{center}
\textit{Note}: A short description of $\lambda_{x_j},
\varepsilon_{x_j}, \delta_{x_j}, \xi_{x_j}$ appears in section 2.1.
For a full
description and examples - see \cite{MoTe2}.\\
$\xi_{x_1} = z_{3,4}\\\varphi_M(\ell(\gamma_1)) = Z^2_{3,4}$\\[1ex]
$\xi_{x_2} = z_{4',5}$ ( $\Delta\!<\!3,4\!>$ does not affect this
path)
$\\\varphi_M(\ell(\gamma_2)) = Z^4_{4',5}$\\[1ex]
$\xi_{x_3} = $\,
\epsfig{file=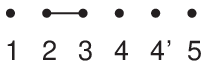}$\xrightarrow[\Delta<3,4>]{\Delta^2<4',5>}$\,
\epsfig{file=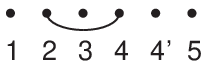}$ =
z_{2\,4}\\[1ex]\varphi_M(\ell(\gamma_3)) = Z^4_{2\,4}$\\[1ex]
$\xi_{x_4} = $\,
\epsfig{file=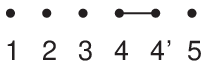}$\xrightarrow[\Delta^2<4',5>]{\Delta^2<2,3>}$\,
\epsfig{file=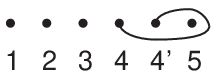}$\xrightarrow{\Delta<3,4>}$\,
\epsfig{file=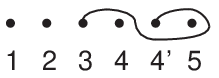}$ =
\underset{\{5\}}{\overset{(4)}{\underline{z}^2_{3\,4'}}}\\[1ex]\varphi_M(\ell(\gamma_4)) =
\underset{\{5\}}{\overset{(4)}{\underline{Z}^2_{3\,4'}}}$\\[1ex]
$\xi_{x_5} = $\,
\epsfig{file=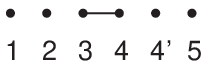}$\xrightarrow[\Delta^2<2,3>]{\Delta<4,4'>}$\,
\epsfig{file=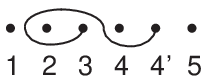}$\xrightarrow{\Delta^2<4',5>}$\,
\epsfig{file=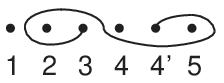}$\xrightarrow{\Delta<3,4>}$\,\\
\epsfig{file=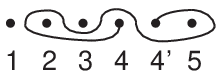}$ =
\hat{\hat{z}}_{4\,4'}\\[1ex]\varphi_M(\ell(\gamma_5)) =
\hat{\hat{Z}}_{4\,4'}$\\[1ex]
$\xi_{x_6}, \xi_{x_{6'}} = $\,
\epsfig{file=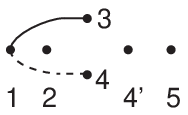}$\xrightarrow[\Delta<4,4'>
\atop \Delta^2<2,3>]{\Delta^{\frac{1}{2}}_{IR}<2>}$\,
\epsfig{file=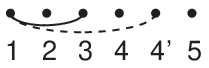}$\xrightarrow[\Delta<3,4>]{\Delta^2<4',5>}$\,
\epsfig{file=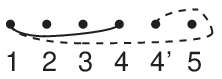}\\[1ex] So - $\xi_{x_6'} =
\underline{z}_{1\,4},\,\xi_{x_6} =
\underset{\{5\}}{\underline{z}_{1\,4'}}$,\\and by Corollary 2.2  -
$\varphi_M(\ell(\gamma_6)\ell(\gamma_{6'})) =
\underset{\{5\}}{\underline{Z}^2_{1\,4'}}\underline{Z}^2_{1\,4}$\\[1ex]
$\xi_{x_7} = $\,
\epsfig{file=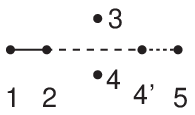}$\xrightarrow{\Delta^2<1,3>}$\,
\epsfig{file=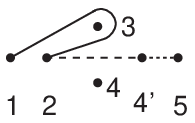}$\xrightarrow{\Delta^{\frac{1}{2}}_{IR}<2>}$\,
\epsfig{file=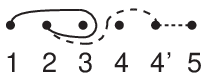}$\xrightarrow{\Delta<4,4'>}$\,
\epsfig{file=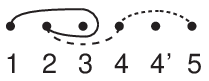}$\xrightarrow{\Delta^2<2,3>}$\,
\epsfig{file=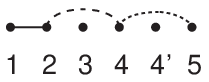}$\xrightarrow{\Delta^2<4',5>}$\,
\epsfig{file=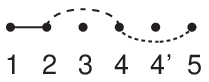}$\xrightarrow{\Delta<3,4>}$\,
\epsfig{file=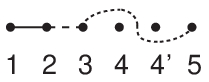}$ =
\Delta\!<\!1,2,3,5\!>^{Z^{-2}_{4,5}}$. Note that in the first
computation of $\xi_{x_7}$ we used Corollary 2.2; therefore, \\[1ex]$\varphi_M(\ell(\gamma_7)) =
\Delta^2\!<\!1,2,3,5\!>^{Z^{-2}_{4,5}}.$\quadf\quadf\quadf\,\,$\square$

 The next
relevant regeneration which affects the neighborhood of $v$ occurs
when we regenerate lines 2 and 3. Note that in a small neighborhood
of $\bigcup\limits_{i=1,2 \atop 3,5}L_i$ the regeneration process
was already treated in \cite{MoTe4}, since the local configuration
of the lines is as in figure 6, and this is exactly the situation
described in \cite{MoTe4}.
\begin{center}
\epsfig{file=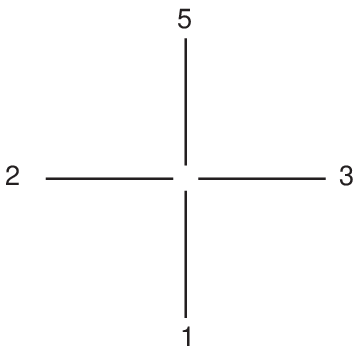}\\
(figure 6)
\end{center}

\begin{corollary}After the described regeneration, the local braid monodromy of
a neighborhood of $v$
is given by-\\
\begin{center}
$\tilde{\tilde{\varphi}} =
Z^2_{3\,4}\,Z^{(3)}_{4',5\,5'}Z^4_{2\,4}\,
\underset{\{5'\}}{\overset{(4)}{\underline{Z}^2_{3\,4'}}}\,\hat{\hat{Z}}_{4\,4'}\,
\underset{\{5'\}}{\underline{Z}^2_{1'\,4'}}\,
\underset{\{5'\}}{\underline{Z}^2_{1\,4'}}\,
\underline{Z}^2_{1'\,4}\, \underline{Z}^2_{1\,4}\,
(B)^{Z^{-2}_{4\,5}Z^{-2}_{4\,5'}}$
\end{center} where
$\hat{\hat{Z}}_{4\,4'} = H(\hat{\hat{z}}_{4\,4'})$, and
$\hat{\hat{z}}_{4\,4'}$ is the path represented by:
\begin{center}
\epsfig{file=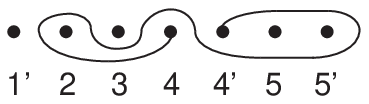}\\
(figure 7)
\end{center}
and $B\,:=\,F(F)_{\rho^{-1}}$ where-\\
$F =
Z^4_{1',2}\,Z^4_{3\,5}\,\tilde{Z}_{2\,3}\,\overset{(2)}{Z^2_{1',5}}\,\underset{(4-4')}{\bar{Z}^2_{1',5'}}\\
\rho = Z_{1\,1'}Z_{5\,5'}$ and $\tilde{Z}_{2,3}$ is represented by
\begin{center}
\epsfig{file=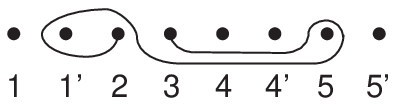}\\
(figure 8)
\end{center}
\end{corollary}

\textbf{Proof:}\, After regenerating $\bigcup\limits_{i=1,2 \atop
3,5}L_i$ in a small neighborhood $U'$ of $v$, $L_2$ and $L_3$ turn
into conics: $Q_2, Q_3$ and $L_1, L_5$ double themselves into
parallel lines $L_1, L_{1'},L_5, L_{5'}$, s.t. $L_1,L_{1'}$ is
tangent to $Q_3$, and $L_5,L_{5'}$ are tangent to $Q_3$. So by
\cite{MoTe4}, Lemma 6, when examining what happens in the process of
the regeneration to $\Delta^2\langle1,2,3,5\rangle$ (in
$\tilde{\varphi}$), its local braid monodromy is $B$. Therefore, in
the expression which represents the local braid monodromy of our
current situation around $v, \Delta^2\langle1,2,3,5\rangle$ is
replaced by $B$. The other changes follow from the regeneration
rules, as $L_1$ and $L_5$ are doubled. Therefore we get
$\tilde{\tilde{\varphi}}$, by
replacing in $\tilde{\varphi}$ the following:\\
\begin{enumerate}
\item $Z^4_{4'\,5}$ into $Z^{(3)}_{4',5\,5'}$ (third regeneration rule)
\item $\underline{Z}^2_{1,4}$ into $\underline{Z}^2_{1',4}\,\underline{Z}^2_{1,4}$ (second regeneration rule)
\item $\underset{\{5\}}{\underline{Z}^2_{1,\,4'}}$ into
$\underset{\{5\}}{\underline{Z}^2_{1',\,4'}}\,\underset{\{5\}}{\underline{Z}^2_{1,\,4'}}$
(second regeneration
rule)\quadf\quad\quad\quad$\square$\end{enumerate}

In the final regeneration that affects the neighborhood of $v$ , the
conics $Q_2,Q_3$ are doubled. Therefore, we have the following
proposition:
\begin{corollary}
The local braid monodromy after the final regeneration around $v$
is given by\\
$$\tilde{\tilde{\tilde{\varphi}}} =
Z^2_{3'\,4}\,Z^2_{3\,4}\,Z^{(3)}_{4',5\,5'}Z^{(3)}_{2\,2',4}\,
\underset{\{5'\}}{\overset{(4)}{\underline{Z}^2_{3'\,4'}}}\,
\underset{\{5'\}}{\overset{(4)}{\underline{Z}^2_{3\,4'}}} \,
\hat{\hat{Z}}_{4\,4'}\,
\underset{\{5'\}}{\underline{Z}^2_{1'\,4'}}\,
\underset{\{5'\}}{\underline{Z}^2_{1\,4'}}\,
\underline{Z}^2_{1'\,4}\, \underline{Z}^2_{1\,4}\,
(\tilde{B})^{\bullet\bullet}$$ where $\hat{\hat{Z}}_{4\,4'}$
corresponds to the path
\begin{center}
\epsfig{file=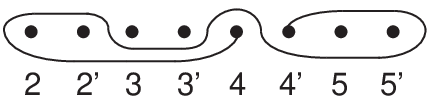}\\
(figure 9)
\end{center}
$(\,)^{\bullet\bullet}$ is conjugation by
$Z^{-2}_{4,5}Z^{-2}_{4,5'}$
\\and $\tilde{B} = \tilde{F}(\tilde{F})_{\rho^{-1}}$ where\\
$\tilde{F} =
Z^{(3)}_{1',2\,2'}\,Z^{(3)}_{3\,3',5}\,\check{Z}_{2'\,3}\,\check{Z}_{2\,3'}\,
\overset{(2-2')}{Z^2_{1',5}}\,\underset{(4-4')}{\bar{Z}^2_{1',5'}}\\
\rho = Z_{1\,1'}Z_{5\,5'}$\\and
$\check{Z}_{2\,3'}\,,\check{Z}_{2'\,3}$ are:
\begin{center}
\epsfig{file=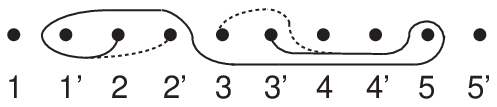}\\
(figure 10)
\end{center}

\end{corollary}
\textbf{Proof:}\, All the changes follow from the regeneration
rules.\\We get $\tilde{\tilde{\tilde{\varphi}}}$ by replacing in
$\tilde{\tilde{\varphi}}$ the following:
\begin{enumerate}
\item $Z^2_{3,4}$ ($\underset{\{5'\}}{\overset{(4)}{\underline{Z}^2_{3\,4'}}}$)
by $Z^2_{3',4}Z^2_{3,4}$ (resp.
$\underset{\{5'\}}{\overset{(4)}{\underline{Z}^2_{3'\,4'}}}
\underset{\{5'\}}{\overset{(4)}{\underline{Z}^2_{3\,4'}}}$) (second
regeneration rule)
\item $Z^4_{2,4}$
by $Z^{(3)}_{2\,2',4}$ (third regeneration rule);
\end{enumerate}
and we get $\tilde{B}\, (\tilde{F})$ by replacing in $B \,($resp.$\,
F)$ the following:
\begin{enumerate}
\item $Z^4_{1',2}$ ($Z^4_{3,5}$) by $Z^{(3)}_{1',2\,2'}$ (resp. $Z^{(3)}_{3\,3',5}$) (third regeneration rule)
\item $\tilde{Z}_{2,3}$ by $\check{Z}_{2'\,3}\,\check{Z}_{2\,3'}$ (first regeneration rule)\quadf\quadf$\square$
\end{enumerate}

\subsection[The second case]{The second case}

The second case of the 5--point regeneration that we deal with is
the braid monodromy factorization that we get from regenerating the
following arrangement of 5 planes corresponding to the angular
sectors of the figure:
\begin{center}
\epsfig{file=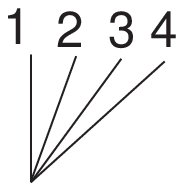}\\
(figure 11)
\end{center}
when first we regenerate line number 4, then line 3, etc. But
instead of looking at a particular case, we examine the general case
where we have a line arrangement of $k$ lines, as in the following
figure:
\begin{center}
\epsfig{file=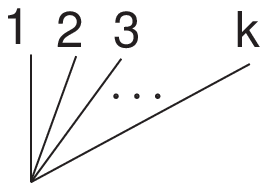}\\
(figure 12)
\end{center}
when first we regenerate line number $k$ (i.e., $\e_k$), then line
$k-1$, etc. We denote $v$ to be  the point of the intersection of
all the lines.

When we regenerate $\e_k$, this line turns into a conic $\tp_{k,k'}$
(by Corollary 2.1) which is tangent to $\e_{k-1}$. In a local
neighborhood of $v$, the real part of this configuration of the
lines $\e_1,...,\e_{k-1}$ and the conic $\tp_{k,k'}$ is as in the
following figure:
\begin{center}
\epsfig{file=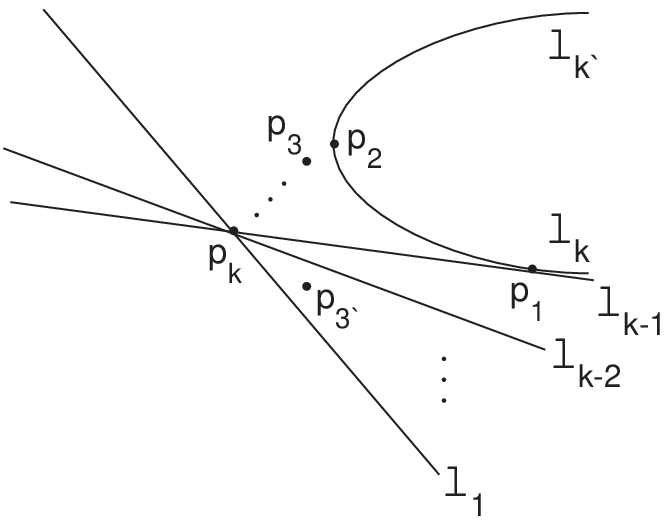}\\
(figure 13)
\end{center}

Note that if we denote the slope of $\e_i$ by $s_i$, then for two
lines - $\e_i, \e_j$ ($1 \leq i < j \leq k-2$ and thus $s_i < s_j$)
we have $\Re(\e_i \cap \tp_{k,k'}) < \Re(\e_j \cap \tp_{k,k'})$.

We shall now compute the braid monodromy factorization of figure 13.
\begin{corollary}
After regenerating $\e_k$, the braid monodromy factorization in a
local neighborhood of $v$ is:
$$B_k = Z^4_{k-1,k}\underset{\{k-1\}}{\bar{Z}_{k,k'}}
\prod\limits_{i=3}^k(Z^2_{k-i+1,k'}Z^2_{k-i+1,k})\Delta^2\langle1,k-1\rangle$$
\end{corollary}
\textbf{Proof:} After the regeneration, $\e_k$ turns into a conic
$\tp_{k,k'}$ and thus, by figure 8, we have the following singular
points with respect to $\p$  (the projection to the X-axis):\\
$p_1$: tangent point of $\e_{k-1}$ and $\tp_{k,k'}$.\\
$p_2$: branch point of $\tp_{k,k'}$.\\
$p_i,p_{i'}: \forall 3 \leq i \leq k$ the (complex) intersection
points (nodes) of $\e_{k-i+1}$ \\ and $\tp_{k,k'}$.\\
$p_{k+1}$: the intersection point of the lines $\e_1,...\e_{k-1}$.\\

We condense the needed data in the following table:\\
\begin{center}
\begin{tabular}{cccc}
$j$ & $\lambda_{x_j}$ & $\varepsilon_{x_j}$ & $\delta_{x_j}$
\\[0.5ex]\hline
1 & $<k-1,k>$ & 4 & $\Delta^2<k-1,k>$\\
2 & $<k,k'>$ & 1 & $\Delta_{IR}<k-1>$\\
3,3' & $<k-2,k>,<k-2,k'>$ & 2 & $\Delta^2<k-2,k>$\\
4,4' & $<k-3,k>,<k-3,k'>$ & 2 & $\Delta^2<k-3,k>$\\
\vdots & $ \vdots $ & \vdots & $ \vdots $\\
k,k' & $<1,k>,<1,k'>$ & 2 & $\Delta^2<1,k>$\\
k+1 & $<1,k-1>$ & 2 & $ - $
\end{tabular}
\end{center}
Therefore\\
$\xi_{x_1} = z_{k-1,k}\\\varphi_M(\ell(\gamma_1)) = Z^4_{k-1,k}$\\[1ex]
$\xi_{x_2} = $\,
\epsfig{file=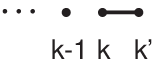}$\xrightarrow{\Delta^2<k-1,k>}$\,
\epsfig{file=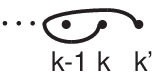}$ =
\underset{\{k-1\}}{\bar{z}_{k,k'}}\\[1ex]\varphi_M(\ell(\gamma_2)) =
\underset{\{k-1\}}{\bar{Z}_{k,k'}}$\\[1ex]
$\xi_{x_3} =
$\,\epsfig{file=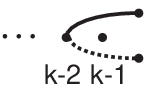}$\xrightarrow{\Delta_{IR}<k-1>}$\,
\epsfig{file=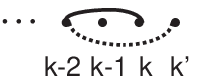}$\xrightarrow{\Delta^2<k-1,k>}$\,
\epsfig{file=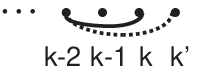}\\[1ex]
$\varphi_M(\ell(\gamma_3)\ell(\gamma_{3'})) = Z^2_{k-2,k'}Z^2_{k-2,k}$\\[1ex]
$\xi_{x_4} =
$\,\epsfig{file=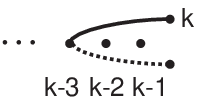}$\xrightarrow{\Delta^2<k-2,k>}$\,
\epsfig{file=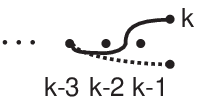}$\xrightarrow{\Delta_{IR}<k-1>}$\,
\epsfig{file=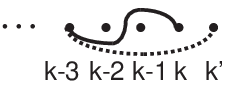}$\xrightarrow{\Delta^2<k-1,k>}$\,
\epsfig{file=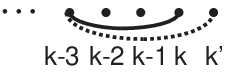}\\[1ex]
$\varphi_M(\ell(\gamma_4)\ell(\gamma_{4'})) = Z^2_{k-3,k'}Z^2_{k-3,k}$.\\[1ex]
 Thus, for $3 \leq i \leq k$:\\
$\xi_{x_i} =
$\,\epsfig{file=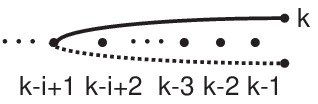}$\xrightarrow{\Delta^2<k-i+2,k>}$\,
\epsfig{file=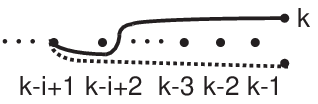} $\cdots
\\[1ex] \xrightarrow{\Delta^2<k-3,k>}$
\epsfig{file=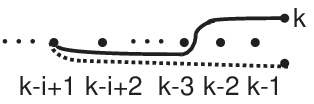}$\xrightarrow{\Delta^2<k-2,k>}$\,
\epsfig{file=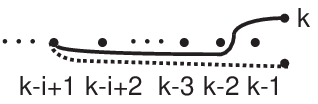}$\xrightarrow{\Delta_{IR}<k-1>}\\[1ex]$\,
\epsfig{file=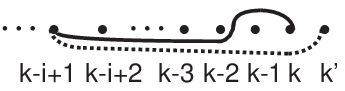}$\xrightarrow{\Delta^2<k-1,k>}$\,
\epsfig{file=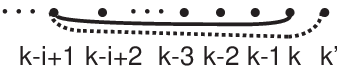}\\[1ex]
$\varphi_M(\ell(\gamma_i)\ell(\gamma_i')) = Z^2_{k-i+1,k'}Z^2_{k-i+1,k}$\\[1ex]
$\xi_{x_{k+1}} =
$\,\epsfig{file=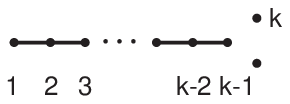}$\xrightarrow{\Delta^2<1,k>}$\,
\epsfig{file=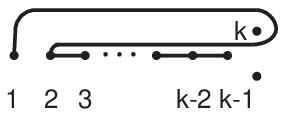}$\xrightarrow{\Delta^2<2,k>}$\,
\epsfig{file=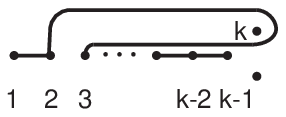}$\cdots \xrightarrow{\Delta^2<k-2,k>}$
\epsfig{file=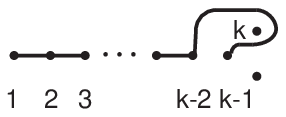}$ \xrightarrow{\Delta_{IR}<k-1>}$
\epsfig{file=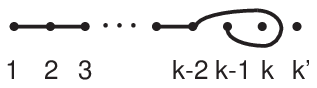}$ \xrightarrow{\Delta^2<k-1,k>}$
\epsfig{file=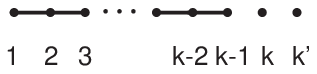}
\\[1ex]$\varphi_M(\ell(\gamma_{k+1})) = \Delta^2\!<\!1,k-1\!>\\
$ Note when computing the braid monodromy factorization in the
neighborhood of the complex points we used Corollary 2.2.
\quadf\quadf\!\!\!$\Box$

The next step is to regenerate $\e_{k-1}$ into a conic:
$\tp_{(k-1),(k-1)'}$. So we have the following:
\begin{corollary}
After the regeneration of $\e_{k-1}$, the braid mondromy
factorization in a local neighborhood of $v$ is:
$$B^{(1)}_k = T_k\prod\limits_{i=3}^k(Z^2_{k-i+1,k'}Z^2_{k-i+1,k})B_{k-1}$$
where  $$T_k =
Z^{(3)}_{((k-1),(k-1)'),k}\underset{\{k-1\}}{\bar{Z}_{k,k'}}.$$
\end{corollary}
\textbf{Proof:} All the changes follow from the regeneration rules.

We get $B^{(1)}_k$ by replacing in $B_k$ the following:
\begin{enumerate}
\item $Z^4_{k-1,k}$ by $Z^{(3)}_{((k-1),(k-1)'),k}$ (by the third regeneration rule)
\item $\Delta^2\!<\!1,k-1\!>$ by $B_{k-1}$ (This is implementation of corollary 2.6,
where we have only $k-1$ lines).\quadf\quadf\quadf\quad$\Box$
\end{enumerate}

The next step is the regeneration of $\e_{k-2}$ (which turns into a
conic $\tp_{(k-2),(k-2)'}$).

So
\begin{corollary}
After the regeneration of $\e_{k-2}$, the braid mondromy
factorization in a local neighborhood of $v$ is:
$$B^{(2)}_k = T_k\prod\limits_{i=3}^3Z^2_{k-i+1,(k-i+1)',k,k'}
\prod\limits_{i=4}^k(Z^2_{k-i+1,k'}Z^2_{k-i+1,k})B^{(1)}_{k-1}$$
where $$Z^2_{\alpha,\alpha',\beta,\beta'} = Z^2_{\alpha',\beta'}
Z^2_{\alpha',\beta}Z^2_{\alpha,\beta'}Z^2_{\alpha,\beta}.$$
\end{corollary}
\textbf{Proof:} We get $B^{(2)}_k$ by replacing in $B^{(1)}_k$ the
following:
\begin{enumerate}
\item $Z^2_{k-2,k'}Z^2_{k-2,k}$ by $Z^2_{k-2,(k-2)',k,k'}$ (by the second regeneration rule)
\item $B_{k-1}$ by $B^{(1)}_{k-1}$ (By implementation of corollary
2.7).\quadf$\Box$
\end{enumerate}

From now on, the braid monodromy factorization after regenerating
$\e_{k-3}$ (and then $\e_{k-4}$ etc.) can be found by a recursive
formula, as we apply the regeneration rules. Thus, the braid
monodromy factorization after regenerating $\e_{k-n}$ is
$$B^{(n)}_k = T_k\prod\limits_{i=3}^{n+1}Z^2_{k-i+1,(k-i+1)',k,k'}
\prod\limits_{i=n+2}^k(Z^2_{k-i+1,k'}Z^2_{k-i+1,k})B^{(n-1)}_{k-1}.$$
Naturally, the process ends when there are no lines to regenerate,
that is, after the regeneration of $\e_1$. For example, we examine
the braid monodromy factorization of the 5--point (when $k=4$):\\

$$\Delta^2\!<\!1,4\!>\, \xrightarrow{\text{regenerating}\,\, \e_4}\, B_4 =
Z^4_{3,4}\underset{\{3\}}{\bar{Z}_{4,4'}}
\prod\limits_{i=3}^4(Z^2_{5-i,4'}Z^2_{5-i,4})\Delta^2\!<\!1,3\!>
$$ $$\xrightarrow{\text{regenerating}\,\, \e_3}$$ $$B_4^{(1)} =
T_4\prod\limits_{i=3}^4(Z^2_{5-i,4'}Z^2_{5-i,4})Z^4_{2,3}
\underset{\{2\}}{\bar{Z}_{3,3'}}Z^2_{1,3'}Z^2_{1,3}\Delta^2\!<\!1,2\!>
$$ $$\xrightarrow{\text{regenerating}\,\, \e_2}$$ $$B_4^{(2)} = T_4Z^2_{2,2',4,4'}Z^2_{1,4'}
Z^2_{1,4}T_3Z^2_{1,3'}Z^2_{1,3}Z^4_{1,2}
\underset{\{1\}}{\bar{Z}_{2,2'}}$$
$$\xrightarrow{\text{regenerating}\,\, \e_1}$$ $$B_4^{(3)} =
T_4Z^2_{2,2',4,4'}Z^2_{1,1',4,4'}T_3Z^2_{1,1',3,3'}T_2$$


 Thus,
$B^{(4)}_5$ is the braid monodromy factorization of the fully
regenerated neighborhood of the 5-point.\\\\
\textbf{Remark}: It is easy to prove (using the recursive formula)
that the braid monodromy factorization of the fully regenerated
neighborhood of the $k+1$-point in figure 12 is:
$$B^{(k-1)}_k = \prod\limits^2_{j=k}\Bigg(T_j\prod\limits^1_{m=j-2}Z^2_{m,m',j,j'}\Bigg)$$
where
$$T_j = Z^{(3)}_{((j-1),(j-1)'),j}\underset{\{j-1\}}{\bar{Z}_{j,j'}}.$$

\section[Using the 5-point regeneration: the Hirzebruch surface $F_{2,(2,2)}$]{Using the 5-point regeneration:\\ the Hirzebruch surface $F_{2,(2,2)}$}

In this section we give an example of using the special braid
monodromy factorization of the 5-point, described in subsection 2.2,
in order to find the global braid monodromy factorization of the
branch curve of a generic projection of $F_{2,(2,2)}$.\\
 \textbf{Remark}: The second case of the regeneration (considered in subsection 2.3,
 for arbitrary $k$) appears when we compute, for example, the global braid
 monodromy factorization of
the branch curve of a generic projection of $F_{k,(a,b)}$ when $k >
2$. Note that finding the global braid monodromy of the Hirzebruch
surface $F_{k,(a,b)},\, \forall k,\,\forall a,b > 1$ can be handled
using only the classical 3- and 6-points, the 5-point studied in
subsection 2.2, and the $k+1$-point studied in subsection 2.3.

\subsection[Braid monodromy of the degenerated curve]{Braid monodromy
of the degenerated curve} The configurations below describe the
projective degeneration of $F_{2,(2,2)} = Z^{(0)}\rightsquigarrow
Z^{(1)} \rightsquigarrow \ldots \rightsquigarrow Z^{(9)}$
\\\\
\epsfig{file=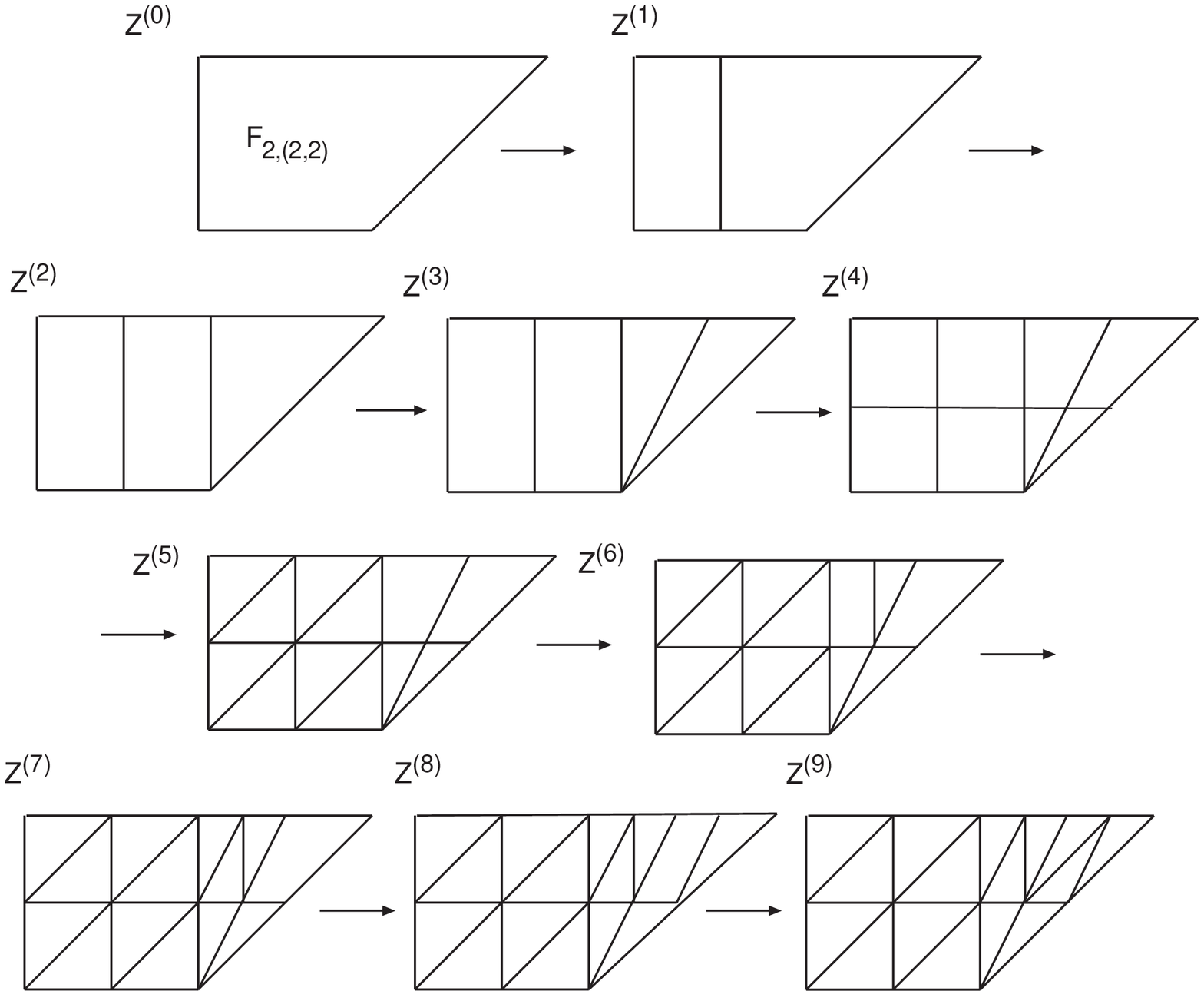,width=12.5cm,height=10cm,angle=0}\\
\centerline{(figure 14)}\\

So $F_{2,(2,2)}$ is degenerated into a union of 16 planes (see
\cite{MoRoTe} for a detailed description of the degeneration of
Hirzebruch surfaces), where the lines represent the intersection of
the planes, and the order of the
vertices is chosen to be lexicographic. See figure 15:\\
\begin{center}
\epsfig{file=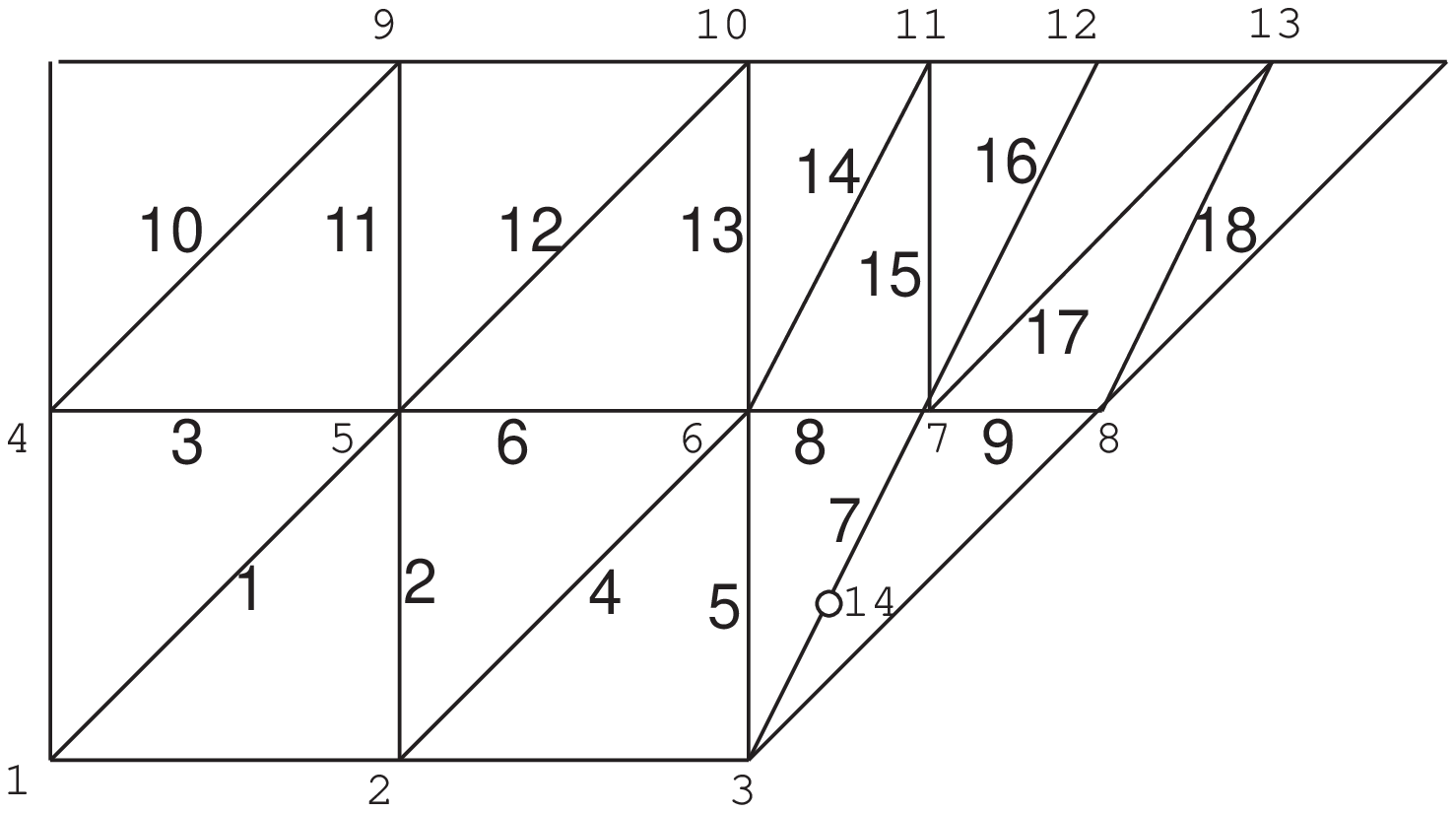,width=7cm,height=4cm,angle=0}\\
(figure 15)
\end{center}
 We denote the numeration of the intersection
lines on $Z^{(9)}$ by $\{\widehat{L_i}\}_{i=1}^{18}$ by the
following method: Let $\widehat{L_1},\widehat{L_2}$ be two edges,
where $\widehat{L_i}$ has vertices $\alpha_i < \beta_i$. Set
$\widehat{L_1} < \widehat{L_2}$ iff $\beta_1<\beta_2$ or
$\beta_1=\beta_2$ and $\alpha_1 < \alpha_2$. Denote also by
$\{\tilde{v}_i\}_{i=1}^{14}$ the intersection points. The appearance
of the point numbered $14$, which is an extra branch point, will be
explained in subsection 3.2.5. Take a generic projection
$\pi^{(i)}:Z^{(i)} \rightarrow \mathbb{CP}^2,\, 0\leq i\leq9$. Let
$S^{(i)}$ be the branch curve in $\mathbb{CP}^2, \varphi^{(i)}$
their braid monodromy, $S^{(i+1)}$ be a degeneration of $S^{(i)}$
(for $0\leq i\leq9$). Let $L_j=\pi^{(9)}(\widehat{L_j}),
j=1\ldots18$. So $S^{(9)}=\bigcup\limits_{j=1}^{18}L_j$; and
$v_j=\pi^{(9)}(\widetilde{v_j}), j=1\ldots14$, so $v_j$ are the
singular points of $S^{(9)}$. Let $C$ be the union of all lines
connecting pairs of the $v_j$-s. $S^{(9)}$ is a subcurve of $C$.
Theorem IX, 2.1, in \cite{MoTe1} gives a full description of the
braid monodromy of $C$: $\Delta^2_C =
\prod\limits_{i=1}^{14}C_i\Delta^2_{v_i}$ with an appropriate
description of $L.V.C$. We use this formula to obtain a description
of $\varphi^{(9)}$ by deleting all factors that involve lines which
do not appear in $S^{(9)}$. Thus, we get $\Delta^2_{S^{(9)}} =
\prod\limits_{i=1}^{14}\widetilde{C_i}\tilde{\Delta}_{v_i}^2$. We
describe each factor separately.\\
$\underline{\tilde{\Delta}}_{v_i}^2$: In $S^{(9)}$, we have 3 points
which are 6-point (points that arise from the intersection of 6
planes), which are $v_5,v_6,v_7$; 8 points which are 3-point, which
are \{$v_j$\}, j=2,3,4,8,9,10,11,13 and two points which are
2-point, which are $v_1,v_{12}$. We denote by $v_{14}$ the extra
branch point, which appears during the regeneration of the line
$L_7$ (see subsection 3.2.5). Since it contributes a factor to the
final braid monodromy factorization, we mention now that the
resulting braid monodromy factorization will be denoted as $\Delta^2
= \prod\limits_{i=1}^{14}C'_i\varphi_i$.

 The local braid monodromies --
$\varphi^{(9)}_j$ ,which are $\tilde{\Delta}_{v_i}^2$, are
introduced and regenerated in the following paragraphs.\\
$\underline{\widetilde{C_i}}$: We get 18 lines in $Z^{(10)}$. Each
line $L_i$ is represented as a pair of its two end vertices. We
define $L_i < L_j$ as above. Define $D_t = \prod\limits_{p<t \atop
L_p\cap L_t=\emptyset}\tilde{Z}^2_{pt}$, where $\tilde{Z}_{pt}$
formulated in \cite{MoTe1} (p. 526). $\tilde{Z}_{pt}^2$ are related
to the parasitic intersections, since they are lines which do not
intersect in $\mathbb{CP}^{15}$ but may intersect in
$\mathbb{CP}^2$. Note that $1\leq p,t \leq 13$, since we do not
include $v_{14}$
in this calculation (see explanation for this in the following passage). Thus: \\
$D_1 = D_2 = D_3 = D_6 = id, D_4 = \bar{Z}^2_{1\,4}\bar{Z}^2_{3\,4},
D_5 = \prod\limits_{p=1}^3\underset{(4)}{\bar{Z}^2_{p\,5}}, D_7
=\prod\limits_{p=1 \atop p\neq5}^6\bar{Z}^2_{p \,7},\\ D_8 =
\prod\limits_{p=1}^3\underset{(7)}{\bar{Z}^2_{p\,8}}, D_9 =
\prod\limits_{p=1}^6\bar{Z}^2_{p \,9} , D_{10} = \prod\limits_{p=1
\atop p\neq3}^9\bar{Z}^2_{p \,10}, D_{11}= \prod\limits_{p=4 \atop
p\neq6}^9\underset{(10)}{\bar{Z}^2_{p\,11}}, D_{12} =
\prod\limits_{p=4 \atop p\neq6}^{10}\bar{Z}^2_{p \,12},\\
D_{13}=\prod\limits_{p=1..3, \atop
7,9..11}\underset{(12)}{\bar{Z}^2_{p,\,13}},
D_{14}=\prod\limits_{p=1..3, \atop 7,9..12}\bar{Z}^2_{p,\,14},
D_{15}= \prod\limits_{p=1..6, \atop
10..13}\underset{(14)}{\bar{Z}^2_{p\,15}}, D_{16}=
\prod\limits_{p=1..6, \atop 10..14}\bar{Z}^2_{p,\,16},\\
D_{17}=\prod\limits_{p=1..6, \atop 10..14}\bar{Z}^2_{p,\,17},
D_{18}=\prod\limits_{p=1..8, \atop
10..16}\underset{(17)}{\bar{Z}^2_{p\,18}},\\$ defining $\tilde{C}_j=
\prod\limits_{V_j\in L_t}D_t$, where $V_j$ is the \emph{small}
vertex among the two vertices of $L_t$, and we get \\
\begin{center}$\tilde{C}_1 = id,\, \tilde{C}_2 = D_4,\, \tilde{C}_3 =
D_5\cdot D_7,\, \tilde{C}_4 = D_{10},\, \tilde{C}_5 = D_{11}\cdot
D_{12},$\\$ \tilde{C}_6=D_8\cdot D_{13}\cdot D_{14},\, \tilde{C}_7 =
D_9\cdot D_{15}\cdot D_{16}\cdot D_{17},\,
\tilde{C}_{8}=D_{18},\,$\\$ \tilde{C}_i = id,$ where $ i =
9,10,\ldots,13$.
\end{center}

As was indicated, the factors $\tilde{C}_j$ correspond to parasitic
intersections. For each point we examine the lines that go through
it, and compute the parasitic intersections with the other lines.
Since we have already looked at the lines passing through $v_3$, we
can ignore the line on which the point $v_{14}$ lies (which is
$L_7$), and by abuse of notation we denote $\tilde{C}_{14} = id$.

\subsection[Local braid monodromy of the regenerated curve]{Local
braid monodromy of the regenerated curve}
\subsubsection[Regeneration of $\tilde{\Delta}^2_{v_j}$]{Computation
and regeneration of $\tilde{\Delta}^2_{v_j}$} We will deal with each
type of point separately.

\subsubsection[3-point]{The 3-point type}

\begin{corollary}
The local braid monodromies $\varphi_2, \varphi_3,
\varphi_4,\varphi_8,\varphi_9,\varphi_{10},\varphi_{11},\varphi_{13}$
are:
\begin{center}
$\varphi_2 = Z^{(3)}_{2\,2',4}\cdot
\tilde{Z}_{4\,4'(2)}\hspace{0.4cm} \varphi_3 =
Z^{(3)}_{5\,5',7}\cdot \tilde{Z}_{7\,7'(5)} $ \\[0.01cm] $ \varphi_4
= Z^{(3)}_{3\,3',10}\cdot \tilde{Z}_{10\,10'(3)}\hspace{0.4cm}
\varphi_8 =
Z^{(3)}_{9\,9',18}\cdot \tilde{Z}_{18\,18'(9)}$ \\[0.01cm] $
\varphi_9 = Z^{(3)}_{10',11\,11'}\cdot
\tilde{Z}_{10\,10'(11)}\hspace{0.4cm} \varphi_{10} =
Z^{(3)}_{12',13\,13'}\cdot
\tilde{Z}_{12\,12'(13)}$ \\[0.01cm] $\varphi_{11} =
Z^{(3)}_{14',15\,15'}\cdot \tilde{Z}_{14\,14'(15)}\hspace{0.4cm}
\varphi_{13} = Z^{(3)}_{17',18\,18'}\cdot \tilde{Z}_{17\,17'(18)},$
\end{center}
where $\tilde{Z}_{i\,i'(j)} = H(\tilde{z}_{i\,i'(j)})$, and \,$
\tilde{z}_{i\,i'(j)}$ is the following path:
\begin{center}
\epsfig{file=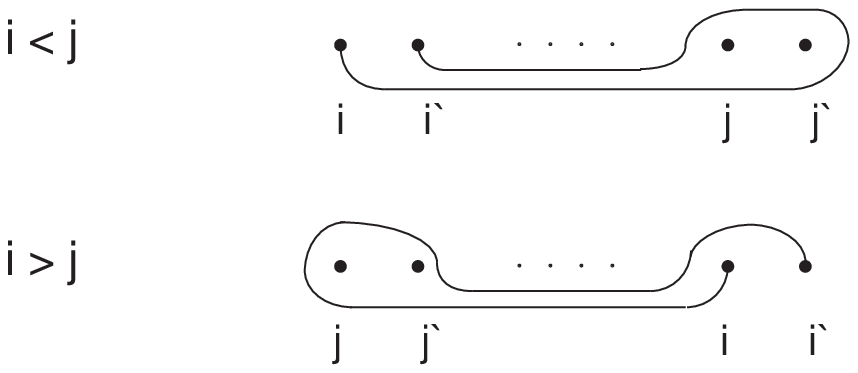,width=6cm,height=3.5cm,angle=0}\\
\emph{(figure 16)}
\end{center}
\end{corollary}
\textbf{Proof:} See \cite{MoTe4}, lemma 1. \\
%
\textbf{Remark:} We will present the representing paths for the
braid monodromy factorization for $\varphi_2$ (figures 17,18) $,
\varphi_9$ (figures 19,20). Note that this configuration of the
paths is the same (with a suitable change of indices) for
$\varphi_3,\varphi_4,\varphi_8$ (resp.
$\varphi_{10},\varphi_{11},\varphi_{13}$):\\\\$Z_{2\,2',4}^{(3)}$:
\begin{center}
\epsfig{file=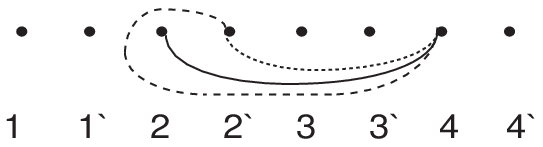}\\
(figure 17)
\end{center}
$\tilde{Z}_{4\,4'(2)}$:
\begin{center}
\epsfig{file=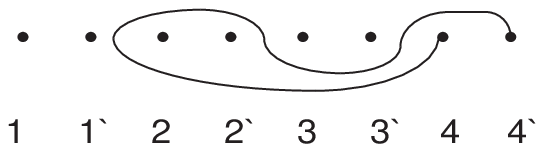}\\
(figure 18)
\end{center}
$Z_{10',11\,11'}^{(3)}$:
\begin{center}
\epsfig{file=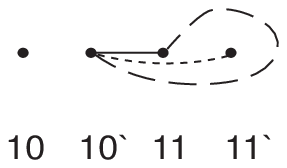}\\
(figure 19)
\end{center}
$\tilde{Z}_{10\,10'(11)}$:
\begin{center}
\epsfig{file=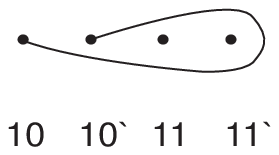}\\
(figure 20)
\end{center}
\subsubsection[6-point]{The 6-point type}
When regenerating $F_{2,(2,2)}$, a new kind of 6-point appears.
Notice that the local numeration of the lines that intersect in
$v_5,v_6$ is as follows:
\begin{center}
\epsfig{file=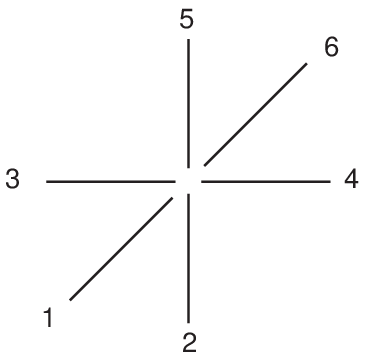}\\
(figure 21)
\end{center}
We will call this kind of 6-point 6-PT1 (6-point type 1). Drawing
(and numerating) the neighborhood of $v_7$ locally, we get:\\\\
\begin{center}
\epsfig{file=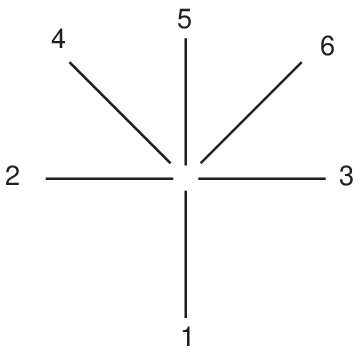}\\
(figure 22)
\end{center}
This kind of 6-point will be called 6-PT2 (6-point type 2). We will
deal first with the regeneration of 6-PT1, which is more familiar,
and then with 6-PT2. \\\\
\textbf{\underline{6-point type 1(6-PT1)}}\\\\ Looking at figure 21,
we see that this configuration of lines was already investigated in
\cite{MoTe4}. Therefore, we cite the main results from there:
\begin{corollary}
The local braid monodromies of $\varphi_5,\varphi_6$ are:\\

$\varphi_5=Z^{(3)}_{1',2\,2'}\tilde{Z}_{12\,12'}Z^{(2)}_{3\,3',12'}(Z^{(2)}_{2\,2',12'})^{\bullet}
\bar{Z}^{(3)}_{6\,6',12}(Z^{(2)}_{3\,3',12})^{\bullet}(Z^{(2)}_{2\,2',12})^{\bullet}
(\hat{F}_{5,1}(\hat{F}_{5,1})_{\rho_5^{-1}})^{\bullet}Z^{(3)}_{11\:11',12}$\\
\begin{center} $\Bigg(\prod\limits_{i=12',12,11'
\atop11,6',6}(Z^2_{1',i})\Bigg)^{\bullet}
\,\bar{Z}^{(3)}_{1',3\,3'}\prod\limits_{i=12',12,11'
\atop11,6',6}(Z^2_{1\,i})\tilde{Z}_{1,1'},$\end{center} where $(
)^{\bullet}$ is the conjugation by the braid induced from the
motion:
\begin{center}
\epsfig{file=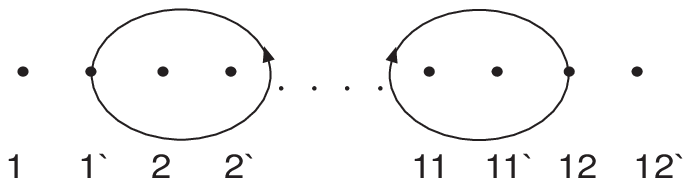}\\
\emph{(figure 23)}
\end{center}
and $\tilde{Z}_{1\,1'}, \tilde{Z}_{12\:12'}$ are
\begin{center}
\epsfig{file=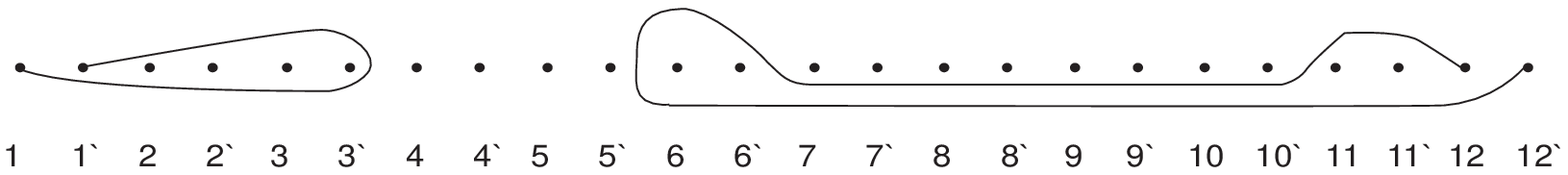}\\
\emph{(figure 24)}
\end{center}
$\rho_5 = Z_{2\,2'}Z_{11\:11'}$\\
$\hat{F}_{5,1}=Z^{(3)}_{2',3\,3'}Z^{(3)}_{6\,6',11}\check{Z}_{3'\,6}\check{Z}_{3\,6'}
\overset{(3-3')}{Z^2_{2',11}}\underset{(4-5', \atop
7-10')}{\bar{Z}^2_{2',11'}}$\\
where $\check{Z}_{3\,6'},\,\check{Z}_{3'\,6}$ are:
\begin{center}
\epsfig{file=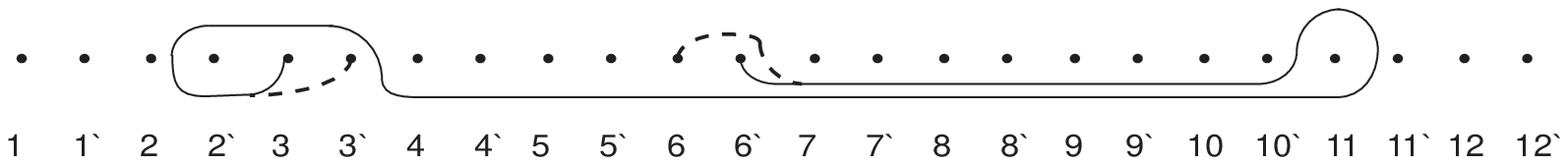}\\
\emph{(figure 25)}
\end{center}
$\varphi_6=Z^{(3)}_{4',5\,5'}\tilde{Z}_{14\,14'}Z^{(2)}_{6\,6',14'}(Z^{(2)}_{5\,5',14'})^{\bullet\bullet}
\bar{Z}^{(3)}_{8\,8',14}(Z^{(2)}_{6\,6',14})^{\bullet\bullet}(Z^{(2)}_{5\,5',14})^{\bullet\bullet}
(\hat{F}_{6,1}(\hat{F}_{6,1})_{\rho_6^{-1}})^{\bullet\bullet}$\\
\begin{center}$Z^{(3)}_{13\:13',14}\Bigg(\prod\limits_{i=14',14,13'
\atop13,8',8}(Z^2_{4',i})\Bigg)^{\bullet\bullet}\,
\bar{Z}^{(3)}_{4',6\,6'}\prod\limits_{i=14',14,13'
\atop13,8',8}(Z^2_{4\,i})\tilde{Z}_{4,4'},$\end{center} where $(
)^{\bullet\bullet}$ is the conjugation by the braid induced from the
motion:
\begin{center}
\epsfig{file=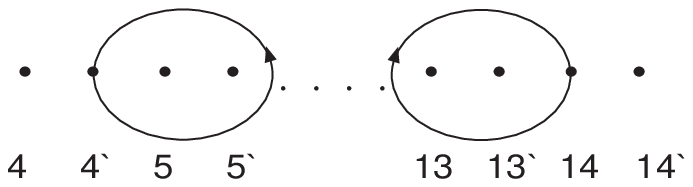}\\
\emph{(figure 26)}
\end{center}
and $\tilde{Z}_{4\,4'}, \tilde{Z}_{14\:14'}$ are
\begin{center}
\epsfig{file=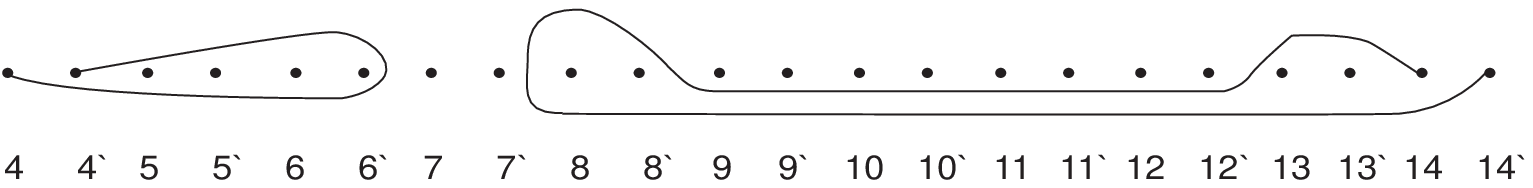}\\
\emph{(figure 27)}
\end{center}
$\rho_6 = Z_{5\,5'}Z_{13\:13'}$\\
$\hat{F}_{6,1}=Z^{(3)}_{5',6\,6'}Z^{(3)}_{8\,8',13}\check{Z}_{6'\,8}\check{Z}_{6\,8'}
\overset{(6-6')}{Z^2_{5',13}}\underset{(7-7', \atop
9-12')}{\bar{Z}^2_{5',13'}}$\\
where $\check{Z}_{6\,8'},\,\check{Z}_{6'\,8}$ are:
\begin{center}
\epsfig{file=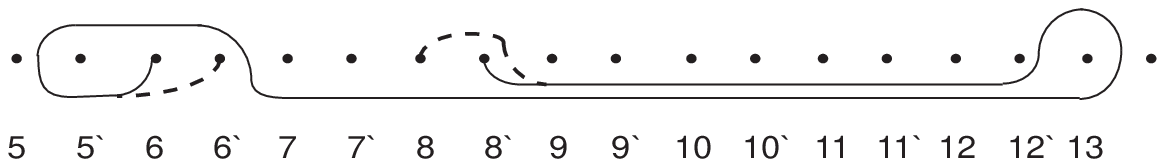}\\
\emph{(figure 28)}
\end{center}
\end{corollary}
\textbf{\underline{6-point type 2 (6-PT2)}}\\\\
We are now dealing with the point $v_7$, that, when numerating the
lines in a local neighborhood of $v_7$, $S^{(9)}$ is as in figure
24. The first regeneration that affects this neighborhood of $v_7$
is the regeneration from $Z^{(9)}$ to $Z^{(8)}$. The line $L_6$, is
regenerated into a conic $Q_6$, that is tangent to $L_3$ and $L_5$.
So, in a small neighborhood of $v_7$, $S^{(8)}$ is as in the
following configuration:
\begin{center}
\epsfig{file=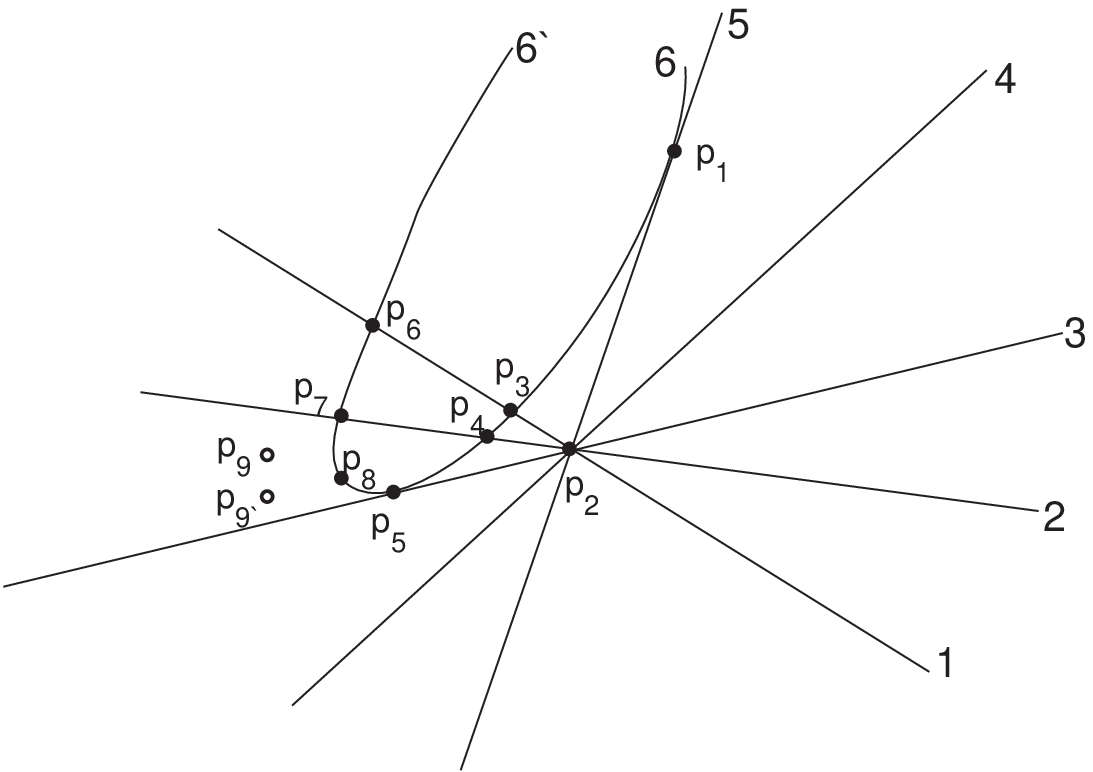}\\
(figure 29)
\end{center}
As we can see, $L_4$ does not intersect the conic in the real part,
so the intersection points of $L_4$ and $Q_6$, $p_9,
p_{9'}\in\mathbb{C}$. By looking at a particular model (where the
conic is $y^2=x, x_{p_2}>0$), it is easy to see that
$\Re{(x_{p_9})}=\Re{(x_{p_{9'}})}<0$. So when calculating the braid
monodromy factorization, we will use Corollary 2.2.

\begin{corollary}
In a neighborhood of $v_7$, the local braid monodromy of $S^{(8)}$
around $v_7$ is given by
\begin{center}
$\varphi_{S^{(8)}} =
Z^4_{5,6}\,(\Delta^2\!<\!1,5\!>)^{Z^2_{5,6}}\,\underset{(5)}{\bar{Z}^2_{1\,6}}
\,\underset{(5)}{\bar{Z}^2_{2\,6}}\,\underset{(5)}{\bar{Z}^4_{3\,6}}\,
\bar{Z}^2_{1\,6'}\,\bar{Z}^2_{2\,6'}\,\tilde{Z}_{6,6'}
\,\tilde{\tilde{Z}}^2_{4,6'}\,\underline{Z}^2_{4,6},$
\end{center}
where the path representing the braid $\tilde{Z}_{6,6'}$ is:
\begin{center}
\epsfig{file=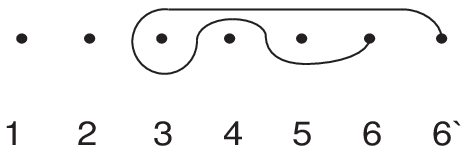}\\
\emph{(figure 30)}
\end{center}
and the path representing the braid $\tilde{\tilde{Z}}^2_{4,6'}$ is
\begin{center}
\epsfig{file=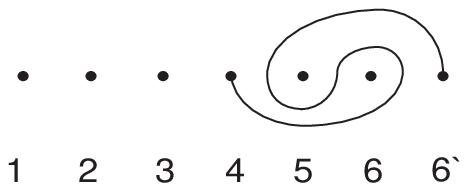}\\
\emph{(figure 31)}
\end{center}
\end{corollary}
\textbf{Proof:} Let $\{p_j\}_{j=1}^9 \cup \{p_{9'}\}$ be the
singular points of a small neighborhood of $v_7$ (see figure 31)
with respect to $\pi_1^{(8)}$ (the projection to the X-axis) as
follows: \\
$p_1, p_6$ - tangent points of $Q_6$.\\ $\{p_3, p_5\},\{p_4,
p_7\},\{p_9, p_{9'}\}$ are the intersection points of $Q_6$ with
$L_1 / L_2 / L_4 $ (resp.).\\ $p_2$ - an intersection point of
$\{L_i\}^5_{i=1}$. \\$p_8$ - the branch point of $Q_6$. \\Let E
(resp. D) be a closed disk on the $X$-axis (resp. $Y$-axis). Let $N
= \{x(p_j) = x_j \:|\: 1\leq j\leq 9 $\:or\:$ j = 9' \}$, s.t.
$N\subset E -
\partial E$. Let $M$ be a real point on the $x$-axis, s.t. $x_j \ll M, \forall x_j \in N,
 1\leq j\leq 9 $\:or\:$ j = 9' \}$. There is a $g$-base
$\ell (\gamma_j)_{j=1}^{9,9'}$ of $\pi_1 (E - N,u)$, s.t. each path
$\gamma_j$ is below the real line and the values of $\varphi_M$
w.r.t this base and $E \times D$ are the ones given in the
proposition. We look for $\varphi_M(\ell(\gamma_j))$ for
$j=1,\ldots,9,9'$. Choose a $g$-base $\ell\{\gamma_j\}_{j=1}^{9,9'}$
as above, and put all the data in the following table:\\
\begin{center}
\begin{tabular}{cccc}
$j$ & $\lambda_{x_j}$ & $\varepsilon_{x_j}$ & $\delta_{x_j}$
\\[0.5ex]\hline
1 & $<5,6>$ & 4 & $\Delta^2<5,6>$\\
2 & $<1,5>$ & 2 & $\Delta<1,5>$\\
3 & $<5,6>$ & 2 & $\Delta<5,6>$\\
4 & $<4,5>$ & 2 & $\Delta<4,5>$\\
5 & $<3,4>$ & 4 & $\Delta^2<3,4>$\\
6 & $<6,6'>$ & 2 & $\Delta<6,6'>$\\
7 & $<5,6>$ & 2 & $\Delta<5,6>$\\
8 & $<4,5>$ & 1 & $\Delta^{\frac{1}{2}}_{IR}<3>$\\
9 & $<2,6>$ & 2 &  - \\
9' & $<2,6'>$ & 2 &  - \\
\end{tabular}
\end{center}
$\xi_{x_1} = z_{5,6}\\\varphi_M(\ell(\gamma_1)) = Z^4_{5,6}$\\[1ex]
$\xi_{x_2} = $\, \epsfig{file=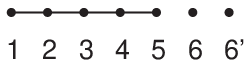}
$\xrightarrow{\Delta^2<5,6>}$\,
\epsfig{file=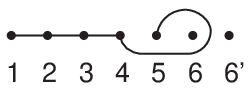}\,
$=\Delta<1,5>^{Z^2_{5,6}}\\[1ex]\varphi_M(\ell(\gamma_2)) =
(\Delta^2<1,5>)^{Z^2_{5,6}}$\\[1ex]
$\xi_{x_3} = $\,
\epsfig{file=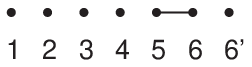}$\xrightarrow{\Delta<1,5>}$\,
\epsfig{file=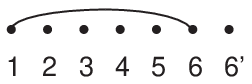}$\xrightarrow{\Delta^2<5,6>}$\,
\epsfig{file=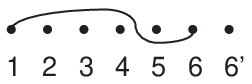}$ =
\underset{(5)}{\bar{z}_{1\,6}}\\\varphi_M(\ell(\gamma_3)) =
\underset{(5)}{\bar{Z}^2_{1\,6}}$\\[1ex]
$\xi_{x_4} = $\,
\epsfig{file=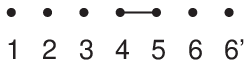}$\xrightarrow{\Delta<5,6>}$\,
\epsfig{file=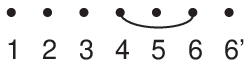}$\xrightarrow{\Delta<1,5>}$\,
\epsfig{file=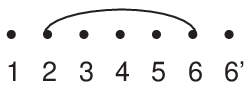}$\xrightarrow{\Delta^2<5,6>}$\,
\epsfig{file=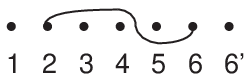}$ =
\underset{(5)}{\bar{z}_{2\,6}}\\[1ex]\varphi_M(\ell(\gamma_4)) =
\underset{(5)}{\bar{Z}^2_{2\,6}}$\\[1ex]
$\xi_{x_5} = $\,
\epsfig{file=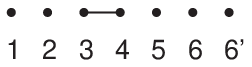}$\xrightarrow{\Delta<4,5>}$\,
\epsfig{file=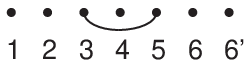}$\xrightarrow{\Delta<5,6>}$\,
\epsfig{file=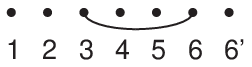}$\xrightarrow{\Delta<1,5>}$\,
\epsfig{file=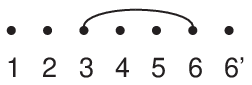}$\xrightarrow{\Delta^2<5,6>}$\,
\epsfig{file=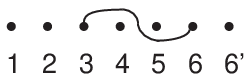}$ =
\underset{(5)}{\bar{z}_{3\,6}}\\[1ex]\varphi_M(\ell(\gamma_5)) =
\underset{(5)}{\bar{Z}^4_{3\,6}}$\\[1ex]
$\xi_{x_6} = $\,
\epsfig{file=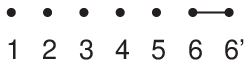}$\xrightarrow[\Delta<4,5>
\atop \Delta<5,6>]{\Delta^2<3,4>}$\,
\epsfig{file=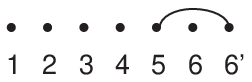}$\xrightarrow{\Delta<1,5>}$\,
\epsfig{file=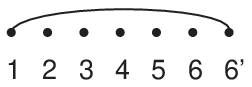}\\[1ex]$ (\Delta^2<5,6> $
does not affect this path) $ =
\bar{z}_{1\,6'}\\\varphi_M(\ell(\gamma_6)) =
\bar{Z}^2_{1\,6'}$\\[1ex]
$\xi_{x_7} = $\,
\epsfig{file=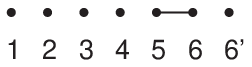}$\xrightarrow{\Delta<6,6'>}$\,
\epsfig{file=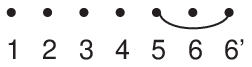}$\xrightarrow[\Delta<4,5>]{\Delta^2<3,4>}$\,
\epsfig{file=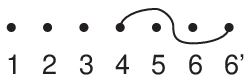}$\xrightarrow{\Delta<5,6>}$\,
\epsfig{file=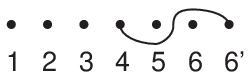}$\xrightarrow{\Delta<1,5>}$\,
\epsfig{file=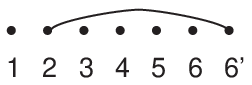}$ =
\bar{z}_{2\,6'}\\[1ex]\varphi_M(\ell(\gamma_7)) =
\bar{Z}^2_{2\,6'}$\\[1ex]
$\xi_{x_8} = $\,
\epsfig{file=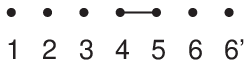}$\xrightarrow[\Delta<5,6>]{\Delta<6,6'>}$\,
\epsfig{file=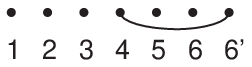}$\xrightarrow{\Delta^2<3,4>}$\,
\epsfig{file=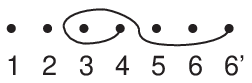}$\xrightarrow{\Delta<4,5>}$\,
\epsfig{file=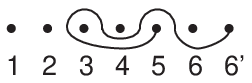}$\xrightarrow{\Delta<5,6>}$\,
\epsfig{file=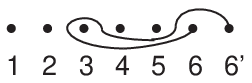}$\xrightarrow{\Delta<1,5>}$\,
\epsfig{file=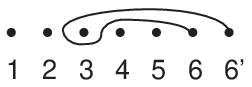}$\xrightarrow{\Delta^2<5,6>}$\,
\epsfig{file=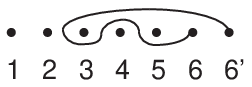}$ = \tilde{z}_{6,6'}
\\[1ex]\varphi_M(\ell(\gamma_8)) = \tilde{Z}_{6,6'}$\\[1ex]
$\xi_{x_9},\,\xi_{x_9'} = $\,
\epsfig{file=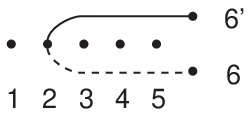}$\xrightarrow{\Delta^{\frac{1}{2}}_{IR}<3>}$\,
\epsfig{file=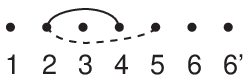}$\xrightarrow[\Delta<6,6'>
\atop \Delta^2<3,4>]{\Delta<5,6>}$\,
\epsfig{file=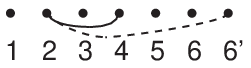}$\xrightarrow[\Delta<5,6>]{\Delta<4,5>}$\,
\epsfig{file=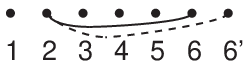}$\xrightarrow{\Delta<1,5>}$\,
\epsfig{file=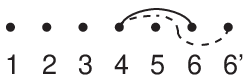}$\xrightarrow{\Delta^2<5,6>}$\,
\epsfig{file=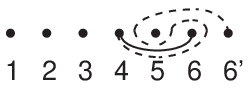} \\
so - $\xi_{x_9'} = \underline{z}_{4,6},\, \xi_{x_9} =
\tilde{\tilde{z}}_{4,6'}$\\and by Corollary 2.2,
$\varphi_M(\ell(\gamma_9)\ell(\gamma_{9'})) =
\tilde{\tilde{Z}}^2_{4,6'}\underline{Z}^2_{4,6}. \quadf\quadf \Box$\\

Note that we are now dealing with a situation described in the
Section 2. That is, a local neighborhood of $v_7$ will look like
figure 1 (in subsection 2.2). Therefore, we can use Corollary 2.5,
which describes what happens to the braid monodromy factorization
after all the regenerations (that is, the local braid monodromy of
$S^{(0)}$).\\

\begin{corollary}
The local braid monodromy of $S^{(0)}$ around $v_7$
is given by,\\
$\varphi_7 = Z^{(3)}_{5\,5',6}\,({\varphi})^{\bullet}\,
\prod\limits_{i=1,\atop 1',2,2'}\underset{(5-5')}{\bar{Z}^2_{i\,6}}
\overset{(4-4')}{\underline{Z}^{(3)}_{3\,3',6}}\, \prod\limits_{i=1,
\atop 1',2,2'}\bar{Z}^2_{i\,6'}\, \tilde{Z}_{6,6'}\,
\tilde{\tilde{Z}}^2_{4',6'} \tilde{\tilde{Z}}^2_{4,6'}\,
\underline{Z}^2_{4',6}\, \underline{Z}^2_{4,6}, \,$\\
where $\tilde{\tilde{Z}}_{4\,6'}, \tilde{\tilde{Z}}_{4'\,6'}$ are as
in the following figure
\begin{center}
\epsfig{file=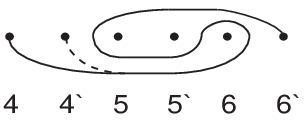}\\
\emph{(figure 32)}
\end{center}
 $\tilde{Z}_{6\,6'}$ corresponds to the
path:
\begin{center}
\epsfig{file=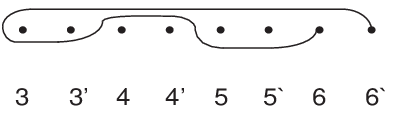}\\
\emph{(figure 33)}
\end{center}
$(\,)^{\bullet}$ is conjugation by $Z^2_{5',6}Z^2_{5,6}$;\\
$\overset{(4-4')}{\underline{Z}^{(3)}_{3\,3',6}}$ is represented by
the 3 following paths:
\begin{center}
\epsfig{file=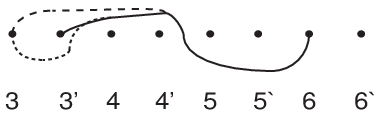}\\
\emph{(figure 34)}
\end{center}
and ${\varphi} =
Z^2_{3'\,4}\,Z^2_{3\,4}\,Z^{(3)}_{4',5\,5'}Z^{(3)}_{2\,2',4}\,
\underset{\{5'\}}{\overset{(4)}{\underline{Z}^2_{3'\,4'}}}\,
\underset{\{5'\}}{\overset{(4)}{\underline{Z}^2_{3\,4'}}} \,
\hat{\hat{Z}}_{4\,4'}\,
\underset{\{5'\}}{\underline{Z}^2_{1'\,4'}}\,
\underset{\{5'\}}{\underline{Z}^2_{1\,4'}}\,
\underline{Z}^2_{1'\,4}\, \underline{Z}^2_{1\,4}\,
(\tilde{B})^{\bullet\bullet},$
\\where $\hat{\hat{Z}}_{4\,4'}$ is as in Corollary 2.5,
$(\,)^{\bullet\bullet}$ is conjugation by
$Z^{-2}_{4,5}Z^{-2}_{4,5'}$
and\\ $\tilde{B} = \tilde{F}(\tilde{F})_{\rho^{-1}}$ where\\
$\tilde{F} =
Z^{(3)}_{1',2\,2'}\,Z^{(3)}_{3\,3',5}\,\check{Z}_{2'\,3}\,\check{Z}_{2\,3'}\,
\overset{(2-2')}{Z^2_{1',5}}\,\underset{(4-4')}{\bar{Z}^2_{1',5'}}\\
\rho = Z_{1\,1'}Z_{5\,5'}$\\and
$\check{Z}_{2\,3'}\,,\check{Z}_{2'\,3}$ are as in Corollary 2.5.
\end{corollary}
\textbf{Proof:} All the changes follow from the regeneration rules
and from Corollary 2.5. Thus, we get $\varphi_7$ by replacing in
$\varphi_{S^{(8)}}$ the following:
\begin{enumerate}
\item By the third regeneration rule: $Z^{4}_{5\,,6} \rightarrow Z^{(3)}_{5\,5',6} ,
\underset{(5)}{\bar{Z}^4_{3\,6}} \rightarrow
\overset{(4-4')}{\underline{Z}^{(3)}_{3\,3',6}}$
\item By the second regeneration rule: $(\,)^{Z^2_{5,6}} \rightarrow
(\,)^{Z^2_{5',6}Z^2_{5,6}}
,\underset{(5)}{\bar{Z}^2_{1\,6}}\underset{(5)}{\bar{Z}^2_{2\,6}}
\rightarrow \prod\limits_{i=1,\atop
1',2,2'}\underset{(5-5')}{\bar{Z}^2_{i\,6}}
,\\\bar{Z}^2_{1\,6'}\bar{Z}^2_{2\,6'} \rightarrow \prod\limits_{i=1,
\atop 1',2,2'}\bar{Z}^2_{i\,6'},\underline{Z}^2_{4,6} \rightarrow
\underline{Z}^2_{4',6} \,\underline{Z}^2_{4,6}
,\tilde{\tilde{Z}}^2_{4,6'}\ \rightarrow
\tilde{\tilde{Z}}^2_{4',6'}\,\tilde{\tilde{Z}}^2_{4,6'}$
\item By Corollary 2.5: $\Delta^2\!<\!1,5\!> \rightarrow \varphi$.\quad\quadf\quadf\quad\quad$\square$
\end{enumerate}

\textit{Note:} $(\,)^{\bullet\bullet},\,(\,)^{\bullet}$ =
conjugation by the braids induced from the motions:
\begin{center}
\epsfig{file=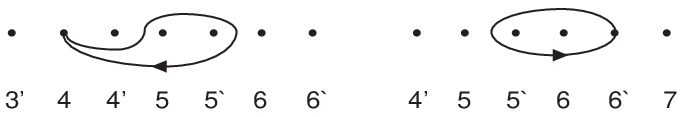}\\
(figure 35)
\end{center}

\textit{Note}: The above computation of $\vp_7$ was done before the
embedding of $B_{12}$ (the braid group with 12 strings, which in all
of the above computation were done) into $B_{36}$ (the braid group
with 36 strings, which in the braid monodromy factorization is
expressed). So we have the following:\\
\textbf{{Remark}}:\,the local braid monodromy of $S^{(0)}$ around
$v_7$, after embedding in $B_{36}$,
is given by\\
$\varphi_7 = Z^{(3)}_{16\,16',17}\,(\tilde{\varphi})^{\bullet}\,
\prod\limits_{i=7,\atop 7',8,8'}\underset{(10-14' \atop
16-16')}{\bar{Z}^2_{i\,17}}
\overset{(15-15')}{\underline{Z}^{(3)}_{9\,9',17}}\,
\prod\limits_{i=7, \atop
7',8,8'}\underset{(10-14')}{\bar{Z}^2_{i\,17'}}
\,\tilde{Z}_{17,17'}\, \tilde{\tilde{Z}}^2_{15',17'}\,
\tilde{\tilde{Z}}^2_{15,17'}\, \underline{Z}^2_{15',17}\,
\underline{Z}^2_{15,17}
$\\
where $\tilde{\tilde{Z}}_{15\,17'},
\tilde{\tilde{Z}}_{15'\,17'},\,\tilde{Z}_{17\,17'}$ are as in the
following figure:
\begin{center}
\epsfig{file=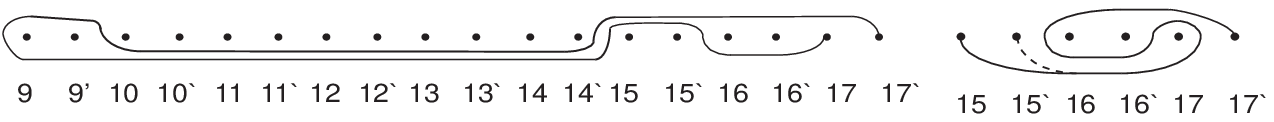}\\
{(figure 36)}
\end{center}
$(\,)^{\bullet}$ is conjugation by $Z^2_{16',17}Z^2_{16,17}$; and \\
$\tilde{\varphi} =
Z^2_{9'\,15}\,Z^2_{9\,15}\,Z^{(3)}_{15',16\,16'}Z^{(3)}_{8\,8',15}\,
\underset{\{16'\}}{\overset{(15)}{\underline{Z}^2_{9'\,15'}}}
\,\underset{\{16'\}}{\overset{(15)}{\underline{Z}^2_{9\,15'}}}\,\hat{\hat{Z}}_{15\,15'}\,
\underset{\{16'\}}{\underline{Z}^2_{7'\,15'}}\,
\underset{\{16'\}}{\underline{Z}^2_{7\,15'}}\,
\underline{Z}^2_{7'\,15}\, \underline{Z}^2_{7\,15}\,
(\tilde{B})^{\bullet\bullet}$ \\where $\hat{\hat{Z}}_{15\,15'}$
corresponds to the path:
\begin{center}
\epsfig{file=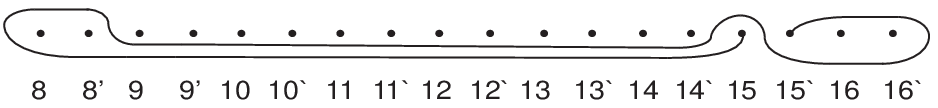}\\
{(figure 37)}
\end{center}
$(\,)^{\bullet\bullet}$ is conjugation by
$Z^{-2}_{15,16}Z^{-2}_{15,16'}$
and\\ $\tilde{B} = \tilde{F}(\tilde{F})_{\rho^{-1}}$ where\\
$\tilde{F} =
Z^{(3)}_{7',8\,8'}\,Z^{(3)}_{9\,9',16}\,\check{Z}_{8'\,9}\,\check{Z}_{8\,9'}\,
\overset{(8-8')}{Z^2_{7',16}}\,\underset{(10-15')}{\bar{Z}^2_{7',16'}}\\
\rho = Z_{7\,7'}Z_{16\,16'}$ and\\
$\check{Z}_{8\,9'}\,,\check{Z}_{8'\,9}$ are:
\begin{center}
\epsfig{file=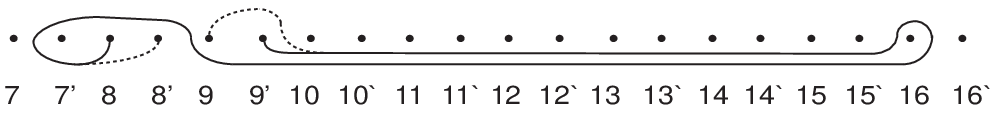}\\
(figure 38)
\end{center}

\subsubsection[Regeneration of $\tilde{C}_i$]{Regeneration of
$\tilde{C}_i$} Performing the regeneration affects also the
parasitic line intersection. Denote by $C'_i$ the braid, which is
created from $\tilde{C}_i$ in the regeneration process. Every
$\tilde{C}_i$ is a product of a 2-degree braid $\bar{Z}^2_{ij}$,
which becomes, as a consequence of the second regeneration rule (see
also \cite{MoTe2}, Proposition 3.2),an 8-degree braid:
$\bar{Z}^2_{i\,i',j\,j'}$ (where we denote by
$\bar{Z}^2_{i\,i',j\,j'} =
Z^2_{i\,j}Z^2_{i'\,j}Z^2_{i\,j'}Z^2_{i'\,j'}$). So \\\\
$C'_i = id ,\:$where $i = 1,3',9,10,...,13\,\,,C'_2 = D'_4 ,\: C'_3
= D'_5D'_7 ,\: C'_4 = D'_{10} ,\\ C'_5 = D'_{11}D'_{12} , C'_6 =
D'_8D'_{13}D'_{14} ,\:
C'_7 = D'_9D'_{15}D'_{16}D'_{17},\:C'_8 = D'_{18}.$\\
where \\\\
$D'_4 = \bar{Z}^2_{1\,1',4\,4'}\bar{Z}^2_{3\,3',4\,4'}, D'_5 =
\prod\limits_{p=1}^3\underset{(4-4')}{\bar{Z}^2_{p\,p',5\,5'}}, D'_7
=\prod\limits_{p=1 \atop p\neq5}^6\bar{Z}^2_{p \,p',7\,7'},\\ D'_8 =
\prod\limits_{p=1}^3\underset{(7-7')}{\bar{Z}^2_{p\,p',8\,8'}}, D'_9
= \prod\limits_{p=1}^6\bar{Z}^2_{p \,p',9\,9'} , D'_{10} =
\prod\limits_{p=1 \atop p\neq3}^9\bar{Z}^2_{p \,p',10\,10'},
D'_{11}= \prod\limits_{p=4 \atop
p\neq6}^9\underset{(10-10')}{\bar{Z}^2_{p\,p',11\,11'}}\\ D'_{12} =
\prod\limits_{p=4 \atop p\neq6}^{10}\bar{Z}^2_{p \,p',12\,12'},
D'_{13}=\prod\limits_{p=1..3, \atop
7,9..11}\underset{(12-12')}{\bar{Z}^2_{p\,p',13\,13'}},
D'_{14}=\prod\limits_{p=1..3, \atop
7,9..12}\bar{Z}^2_{p\,p',14\,14'}, D'_{15}= \prod\limits_{p=1..6,
\atop 10..13}\underset{(14-14')}{\bar{Z}^2_{p\,p',15\,15'}},\\
D'_{16}= \prod\limits_{p=1..6, \atop
10..14}\bar{Z}^2_{p\,p',16\,16'}, D'_{17}=\prod\limits_{p=1..6,
\atop 10..14}\bar{Z}^2_{p\,p',17\,17'},
D'_{18}=\prod\limits_{p=1..8, \atop
10..16}\underset{(17-17')}{\bar{Z}^2_{p\,p',18\,18'}}\\$
\newpage
\subsubsection[2-point]{The 2-point type and extra branch points}

There are two points $v_1,v_{12}$ which are 2-point. The
regeneration around the point $v_1$ yields a conic $Q_1$ (resp.
$Q_{16}$); that is because the pair of the two planes, that on their
intersection the point $v_1$ (resp., and looking locally on the
neighborhood of $v_{12}$) lies, is created in the degeneration
process from a projective nonsingular quadric. But, in order that
these 2--points will contribute a factor to the braid monodromy
factorization, they should satisfy the condition of "Extra branch
points", as described in \cite{Robb}, Section 4. However, only $v_1$
satisfies this condition. Thus $v_{12}$ does not contribute a factor
to this factorization (note that one can check this condition only
after calculating the local factorization around $v_7$). By Artin's
theorem, the induced braid monodromy in the neighborhood of $v_1$ is
$Z_{1 \,1'}$, namely, a counterclockwise halftwist of 1,1' . Thus,
the local braid monodromy  around $v_1$ is $\varphi_1 = Z_{1,\,1'}$.

However, there is another point $v_{14}$ (which is not a 2-point)
that is an extra branch point which is created from the regeneration
of the line $L_7$. This extra branch point contributes to our
factorization the factor $Z_{7,7'}$.

We will now prove the mentioned statements, following the ideas
presented in \cite{Robb}; we also follow its method of numeration,
such that the numeration of the extra branch point is done last. Let
us denote by $\tilde{\Delta} =
\prod\limits_{i=1}^{13}C'_i\varphi_i$. If $\tilde{\Delta}$ is a
braid monodromy factorization, then $\tilde{\Delta} = \Delta^2$ and
$deg(\tilde{\Delta}) = 36\cdot(36-1) = 1260$. However,
$deg(\tilde{\Delta})=1259$ (when assuming that $\varphi_1 =
Z_{1,1'}\,,\varphi_{12}= id$ as will be proved in the following
lemma; the explicit calculation is made in subsection 3.3). This,
there is a missing factor of degree 1.

Define the forgetting homomorphisms: $$f_i:
B_{36}[D,\{1,1',...,18,18'\}] \rightarrow B_2[D,\{i,i'\}],\,\, 1
\leq i\leq 18.$$ It is clear that $\forall i, deg(f_i(\Delta^2))=2$.

\begin{lemma}
\emph{(1)} There exists an extra branch point $v_{14}$ that
contributes a factor $\varphi_{14} = Z_{7,7'}$ to the braid
monodromy
factorization.\\
\emph{(2)} The regeneration of neighborhood of the point $v_{12}$
does not contribute a factor to braid monodromy
factorization.\\
\emph{{(3)}} The regeneration of neighborhood of the point $v_{1}$
does contribute a factor $\varphi_1 = Z_{1,\,1'}$ to braid monodromy
factorization.
\end{lemma}
\textbf{Proof:}\: (1) By Lemma 3.3.3 (or Proposition 3.3.4) in
\cite{RobbT} (see also \cite{Robb}, Section 4), it is enough to
prove that $deg(f_7(\tilde{\Delta}))=1$. The braids coming from the
parasitic intersections are sent by $f_7$ (and by any $f_i$, in
fact) to $id$, so it is enough to look only at the factors
$\varphi_j,\,j \in \{1,...,13\}, \,j \neq 1,12$ that involve braids
 one of whose end points are 7 or $7'$ (we omit the factors
$\varphi_1,\varphi_{12}$ since we do not know yet what are they. In
any case, they do not affect the result, since the do not fulfill
the condition mentioned). The only suitable $j$'s are $j=3$ and
$j=7$. Since $v_3$ is a 3-point, by Lemma 2, \cite{MoTe4},
$deg(f_7(\varphi_3))=1$. Examining $\varphi_7$, it is easy to see
that there exists no braid in the factorization of $\varphi_7$ that
is not sent to $id$ by $f_7$ ; Thus $deg(f_7(\varphi_7))=0$.
Therefore
$deg(f_7(\tilde{\Delta}))=1$.\\\\
(2) In order that the regeneration of the point $v_{12}$ will not
contribute a factor to the braid monodromy factorization, we note
that $v_{12}$ is only on the line $L_{16}$. So we actually have to
check if $\sum\limits_{i=1 \atop i \neq
12,1}^{13}deg(f_{16}(\varphi_i))=2$ to prove the lemma (again, the
braids coming from the parasitic intersections are sent to $id$ by
$f_{16}$). Examining $\varphi_7$, and using Lemma 2,\cite{MoTe4}, we
get that $deg(f_{16}(\varphi_i))=2$ (since
$deg(f_{16}(Z^{(3)}_{16\,16',17})) =
deg(f_{16}((Z^{(3)}_{15',16\,16'}))^{Z^2_{16',17}Z^2_{16,17}}) =
1$), where \\ $\forall 1 \leq i \leq 13, i\neq 1,7,12,\,
deg(f_{16}(\varphi_i)) = 0$.\\
(3) This is done using the same method as in (2), by confirming
that\\
$\sum\limits_{i=2}^{13}deg(f_{1}(\varphi_i))=1.\quadf\quadf\quadf\,\,\quadf\square$

\subsection[Global braid monodromy of the regenerated curve] {Global
braid monodromy of the regenerated curve}

\begin{corollary}
$\Delta^2_{36} = \prod\limits_{i=1}^{14}C'_i\tilde{\varphi_i}$ is a
braid monodromy factorization for $S^{(0)}$, where
$\tilde{\varphi_i} = (\varphi_i)_{h_i}$ for some $h_i\in \langle
Z_{jj'}\, |\,v_i \in L_j\rangle$.
\end{corollary}
\textbf{Proof:}\:Using Proposition VI.2.1 from \cite{MoTe1} on
$S^{(0)}$, we get that $\Delta^2_{36} =
\prod\limits_{i=1}^{14}C'_i\tilde{\varphi_i}\prod b_i$, for some
$h_i\in \langle Z_{jj'}\, |\,v_i \in L_j\rangle$ determined by the
regeneration of the embedding $B_k \hookrightarrow B_{18}$ to
$B_{2k}\hookrightarrow B_{36}$, where $k=1$ when $i=1,14,\,k=2$ when
$i=2,3,4,8,9,10,11,13,\,$ and $k=6$ if $k=5,6,7$ (see the definition
of regeneration of an embedding in \cite{MoTe4}, section 1). And
where $b_i$ are factors corresponding to singularities that are not
covered by $\prod C'_i\tilde{\varphi_i}$, and each $b_i$ is of the
form $Y_i^{t_i}, Y_i$ is a positive halftwist, $0 \leq t_i \leq 3$.
Note that $deg(\tilde{\varphi_i}) = deg(\varphi_i)$. \\We shall
compute $deg\Bigg(\prod\limits_{i=1}^{14}C'_i\varphi_i\Bigg) =
\sum\limits^{14}(deg(C'_i)+deg(\varphi_i))$. So,
$\sum\limits^{14}deg(C'_i) = 800$ (there are 100 factors; each
factor - $Z^2_{i\,i',j\,j'}$ - has degree 8).

 For
$i=2,3,4,8,9,10,11,13 - v_i$ are 3-point (in $S^{(9)}$); if $v_i$ is
a 3-point, then $deg(\varphi_i) = 10$ (by corollary 3.1); then
$\sum\limits_{i; v_i is  \atop  3-point}deg(\varphi_i) = 8\cdot 10 =
80$.

For $i=1,\,v_i$ is 2-point; the line $L_7$ has an extra branch point
(that is - in these cases we have a contribution of a factor of the
form $Z_{i,i'}$), so by the description in section
2.1, these factors contribute to the sum of degrees the addend  2.

For $i=5,6,7,\,v_i$ is 6-PT1 (for $i=5,6$) or a 6-PT2 (for
$i=7$). In any case, $\varphi_i$ includes:\\
\quad 6 factors with degree 1 $\Rightarrow \text{degree} = 6$,\\
\quad 24 factors with degree 2 $\Rightarrow \text{degree} = 48$,\\
\quad 24 factors with degree 3 $\Rightarrow \text{degree} = 72$.\\
So, for a 6-point $v_i,\, deg(\varphi_i)=126$ and $\sum\limits_{i;
v_i is  \atop 6-point}deg(\varphi_i)= 3\cdot 126 = 378$. \\
Therefore - $deg(\prod\limits_{i=1}^{14}C'_i\varphi_i) = 800 + 80 +
2 + 378 = 1260$. Since the degree of
$deg(\prod\limits_{i=1}^{14}C'_i\varphi_i)$ is $36 \cdot 35 = 1260 =
\Delta^2_{36}$ we have $deg(\prod b_i) = 1$, since $\forall i, b_i$
is a positive power of a positive halftwist , we get $b_i = 1\,
\forall i$. So we have $\Delta^2_{36} =
\prod\limits_{i=1}^{14}C'_i\tilde{\varphi_i}$.

\subsubsection[Invariance rules]{Invariance rules}

The aim of this subsection is to prove that indeed $\Delta^2_{36} =
\prod\limits_{i=1}^{14}C'_i{\varphi_i}$. For this, we need to define
a few definitions. We start by defining a Hurwitz move on
$G\times\dots\times G$ ($G$ is a group) or on a set of
factorizations.

 {\textit{Definition}:\, \underline{Hurwitz moves}}:

Let $\vec t= (t_1,\ldots ,t_m)\in G^m$\,. We say that $\vec s
=(s_1,\ldots ,s_m)\in G^m$ is obtained from $\vec t$ by the Hurwitz
move $R_k$ (or $\vec t$ is obtained from $\vec s$ by the Hurwitz
move $R^{-1}_k$) if
$$
s_i = t_i \quad\text{for}\  i\ne k\,,\, k+1\,,\\
s_k = t_kt_{k+1}t^{-1}_k\,,\\
s_{k+1} =t_k\,.
$$
\newpage
\textit{Definition}:\, \underline{Hurwitz move on a factorization}

Let $G$ be a group $t\in G.$  Let  $t=t_1\cdot\ldots\cdot t_m=
s_1\cdot\ldots\cdot s_m$ be two factorized expressions of $t.$ We
say that $s_1\cdot\ldots\cdot s_m$ is obtained from
$t_1\cdot\ldots\cdot t_m$ by a Hurwitz move $R_k$ if $(s_1,\ldots
,s_m)$ is obtained from
$(t_1,\ldots ,t_m)$ by a Hurwitz move $R_k$\,.\\
\textit{Definition}:\, \underbar{Hurwitz equivalence of
factorization}

 Two factorizations are Hurwitz equivalent if they are
obtained  from
each other by a finite sequence of Hurwitz moves.\\
\textit{Definition}:\, \underbar{A factorized expression invariant
under $h$}

 Let $t=t_1\cdot\ldots\cdot t_m$ be a factorized
expression in a group $G$. We say that $t$ is invariant under $h \in
G$ if $(t_1)_h\cdot\ldots\cdot (t_m)_h$ is a Hurwitz equivalent to
$t_1\cdot\ldots\cdot t_m$.

We cite now two lemmas that we will need below.
\begin{lemma}\emph{(see \cite{MoTe1})} If a braid monodromy factorization of
$\Delta^2_{36} = \prod \varphi(\Gamma_i)$ (where $\Gamma_i$ is a
g-base: a free base of $\p(\C^1 - N,u)$ with certain properties; see
\emph{ \cite{MoTe1}} for definition) is invariant under $h$, then
the equivalent factorization $\prod(\varphi(\Gamma_i))_h =: \prod
Z_i$ is also a braid monodromy factorization. That is, $\exists\,$ a
g-base {$\Gamma'_i$} of $\p(\C^1 - N,u)$ s.t. $Z_i =
\varphi(\Gamma'_i)$.
\end{lemma}

\begin{lemma}[Chakiri's Lemma] Let $t=t_1\cdot\ldots\cdot t_m$ be a
factorized expression in a group $G$. Then $t_1\cdot\ldots\cdot t_m$
is invariant under $t^k,\,\forall k\in \mathbb{Z}$.
\end{lemma}

We now look at all the invariance relations that are related to any
kind of point.
\begin{lemma}
$\forall i, \forall (m_j)_{1\leq j \leq 18} \in \Z,\, C'_i$ is
invariant under $\varepsilon =
\prod\limits_{i=1}^{18}Z_{j\,j'}^{m_j}$.
\end{lemma}
\begin{lemma}
$\forall i$, s.t. $v_i$ is a $3$-point, $v_i = L_\alpha \cap L_\beta
\, ,\varphi_i$ is invariant under $\varepsilon$.
\end{lemma}
\textbf{Proof:}\,\cite{MoTe4}, \,Corollary 14.\\
\begin{lemma}
$\forall i$, s.t. $v_i$ is a $2$-point or an extra branch point,
$\varphi_i$ is invariant under $\varepsilon$.
\end{lemma}
\textbf{Proof:}\, the 2-points are $v_1, v_{12}$; consider $v_1$
(recall that $v_{12}$ does not contribute a factor to our
factorization). We have to check if $Z_{1,\,1'}$ is invariant under
$Z_{j,\,j'}\,\forall j$. For $j \neq 1,\,Z_{1,\,1'}$ is invariant
under $Z_{j,\,j'}$ (since the paths corresponding to the braids are
disjoint; thus the braids commute). For $j=1$ we have invariance by
Chakiri's Lemma. For $v_{14}$, we apply the same procedure.\,$\square$\\
\begin{corollary} $\varphi_5$ is invariant under
$(Z_{1,\,1'}Z_{12,\,12'})^{p_1}(Z_{2,\,2'}Z_{11,\,11'})^{q_1}(Z_{3,\,3'}Z_{6,\,6'})^{r_1}$,\\
$\varphi_6$ is invariant under
$(Z_{4,\,4'}Z_{14,\,14'})^{p_2}(Z_{5,\,5'}Z_{13,\,13'})^{q_2}(Z_{6,\,6'}Z_{8,\,8'})^{r_2}$,
\linebreak $\forall p_i, q_i, r_i \in \Z, \,i=1,2$.
\end{corollary}
\textbf{Proof:}\, This is the same as Lemma 15 in \cite{MoTe4}.\\\\
We shall prove now the invariance property for the 6-PT2.
\begin{lemma}
$\varphi_7 = \prod\limits_{i = 8,9 \atop
15,17}Z^{-1}_{i,\,i'}Z^{-2}_{7,\,7'}\Delta^2_{12}$.
\end{lemma}
\textbf{Proof:}\, Let $$L = \{ i,i' \, |\, 1 \leq i \leq 18, i \neq
7,8,9,15,16,17\}$$
$$ G = \{b \in B_{36}[D,\{i,i'\, |\, i=1,\ldots,18\}]\,|\,b\{L\} = \{L\}\}.$$
Denote $\nu : G \ri B_{12}[D,\{i,i'\, |\, i=7,8,9,15,16,17\}]$ the
forgetting homomorphism. Thus - $\nu(\Delta^2_{36}) =
\Delta^2_{12}$; each factor in $C'_i$ contains one of the indices
in $L$, so $\nu(C'_i) = 1\,\forall i$.\\
For $i = 1,2,4,5,9,10,11$ all of the indices in $\tilde{\varphi}_i$
are in $L$, and so $\nu(\tilde{\varphi}_i)=1$ for $i = 1,3,8,11,13\:
v_i$ is a 3-point, when only one index of $\tilde{\varphi}_i$ is in
$L$. So we have (by \cite{MoTe4}, lemma 2) $$\nu(\tilde{\varphi}_3)
= Z_{7,\,7'},\: \nu(\tilde{\varphi}_8) =
Z_{9,\,9'},\:\nu(\tilde{\varphi}_{11}) =
Z_{15,\,15'},\:\nu(\tilde{\varphi}_{13}) = Z_{17,\,17'}$$. For $i =
12, v_i$ is a 2-point that does not contribute a factor to the
factorization. For $i=14, v_i$ is an extra branch point, so
$\nu(\tilde{\varphi}_{14}) = Z_{7,\,7'}$.
$v_6$ is a 6-PT1. All of the factors outside $\hat{F}_{6,1}$ contain
indices in L. So $\nu(\tilde{\varphi}_6) =
\nu(\hat{F}_{6,1}(\hat{F}_{6,1})_{\rho_6^{-1}}) = Z_{8,\,8'}$. Thus
 $\Delta^2_{12} = \nu(\tilde{\varphi}_7)\prod\limits_{i = 8,9 \atop
15,17}Z_{i,\,i'}\,Z^2_{7,\,7'}$, and so $h_7$ commutes with
$\Delta^2_{12}\prod\limits_{i = 8,9 \atop
15,17}Z^{-1}_{i,\,i'}\,Z^{-2}_{7,\,7'}$; Therefore, $
\nu(\tilde{\varphi}_7) = \nu(\varphi_7)$; and since the indices in
$L$ do not appear in $\varphi_7,\:\, \nu(\varphi_7) = \varphi_7$.
Thus,\\ $\varphi_7 = \prod\limits_{i = 8,9 \atop
15,17}Z^{-1}_{i,\,i'}\,Z^{-2}_{7,\,7'}\Delta^2_{12}.\!\!\quadf\quadf\quadf\quadf\square$\\
\begin{corollary}
$\varphi_7$ is invariant under
$I_1(p) = (Z_{8,\,8'}Z_{9,\,9'})^{p}(Z_{15,\,15'}Z_{17,\,17'}Z^2_{7,\,7'})^{p}$\\
$\forall p \in \Z$.
\end{corollary}
\textbf{Proof:} We know that $\varphi_7 = \prod\limits_{i = 8,9
\atop 15,17}Z^{-1}_{i,\,i'}\,Z^{-2}_{7,\,7'}\Delta^2_{12}$. By
Chakiri's Lemma, $\varphi_7$ is invariant under
$\Bigg(\prod\limits_{i = 8,9 \atop
15,17}Z^{-1}_{i,\,i'}\,Z^{-2}_{7,\,7'}\Delta^2_{12}\Bigg)^{-p}$.
Since $\Delta^2_{12}$ is a central element, $\varphi_7$ is invariant
under
$(Z_{8,8'}Z_{9,9'})^p(Z_{15,15'}Z_{17,17'}Z^2_{7,\,7'})^p.\,\quad\quadf\quadf \square$\\

Denote - $\vartheta =
Z^{-2}_{15,16}Z^{-2}_{15,16'}Z^2_{16',17}Z^2_{16,17}$.

\begin{corollary}
$\varphi_7$ is invariant under $I_2(p) = \rho^p =
(Z_{7,7'}Z_{16,16'})^p$
\end{corollary}
\textbf{Proof:}\, Looking at the factors outside
$(\tilde{F}(\tilde{F})_{\rho^{-1}})^\vartheta$, we can see, using
the invariance rules 2 and 3 (\cite{MoTe4}), that outside
$(\tilde{F}(\tilde{F})_{\rho^{-1}})^\vartheta$, the factorization is
indeed invariant under $Z_{7,7'}$ and $Z_{16,16'}$, and by the
invariance remark (v) \cite{MoTe4}, is invariant under $\rho$. So it
is enough to check that
$(\tilde{F}(\tilde{F})_{\rho^{-1}})^\vartheta$ is invariant under
$\rho$. But this is proven exactly in the same way as in Lemma 15,
case 2.2,\cite{MoTe4} (since in our case the point $v_7$ and in the
standard case of the 6--point at \cite{MoTe4}, the
regeneration of the 4--point are the same).$\!\!\quadf\quadf\quadf\quadf\quadf\square$\\

\begin{corollary}
$\varphi_7$ is invariant under $I_3(p) = (Z_{8,8'}Z_{9,9'})^p$.
\end{corollary}
\textbf{Proof:}\, We use the invariance rules (\cite{MoTe4}) when
passing on all the factors of $\varphi_7$. We use invariance rules 2
and 3 when passing on the factors outside
$(\tilde{F}(\tilde{F})_{\rho^{-1}})^\vartheta$ (for example, by
invariance rule number 3,
$(Z^{(3)}_{8,8',15})^{Z^2_{16',17}Z^2_{16,17}}$ is invariant under
$Z_{8,8'}$; by invariance rule number 2, $Z^2_{9',15}Z^2_{9,15}$ is
invariant under $Z_{9,9'}$). When passing on the factors of
$(\tilde{F}(\tilde{F})_{\rho^{-1}})^\vartheta$, we use invariance
rule 3 (for factors of the form $Z^{(3)}_{...}$) and invariance rule
1 (for the factors $Z_{8',9}Z_{8,9'}$ and
$(Z_{8',9}Z_{8,9'})_{\rho^{-1}}$).$\quadf\quadf\quadf\quadf\square$\\

\begin{corollary}
$\varphi_7$ is invariant under
$(Z_{8,\,8'}Z_{9,\,9'})^{q}(Z_{15,\,15'}Z_{17,\,17'}Z^2_{7,\,7'})^{p}(Z_{7,7'}Z_{16,16'})^r,
\\\forall p,q,r \in \Z$.
\end{corollary}
\textbf{Proof:}\, By invariant remark (v) (\cite{MoTe4}),
$\varphi_7$ is invariant under \\$I_3(q-p)\cdot I_1(p)\cdot I_2(r)$,
which is the
desired expression.$\quadf\quad\quad\,\,\,\,\square$\\



\textbf{The Main Result}: As a consequence of the invariance rules,
we can apply them as in \cite{MoTe4}  which remains the same both in
$v_7$ and in the standard 6-point) and get that $\varepsilon(36) :=
\prod\limits_{i=1}^{14}C'_i{\varphi_i}$ is also a braid monodromy
factorization.

Note, that although that the invariance rules for $v_7$ are
different from the invariance rules of the standard 6-point, what
matters , as can be seen in \cite{MoTe4}, Section 4, is that the
invariance rule regarding the horizontal lines in the 6-point (the
two lines that are regenerated last).

\textsc{Michael Friedman, Department of Mathematics,
 Bar-Ilan University, 52900 Ramat Gan, Israel\\}
 \textsl{email}: fridmam@mail.biu.ac.il \\\\
\textsc{Mina Teicher, Department of Mathematics,
 Bar-Ilan University, 52900 Ramat Gan, Israel\\}
 \textsl{email}: teicher@macs.biu.ac.il
\end{document}